\documentclass[a4paper,11pt]{amsart}
\usepackage{amssymb,amsthm,amscd,amsxtra,eucal,fullpage}

\theoremstyle{plain}
\newtheorem{thm}{Theorem}[section]
\newtheorem*{Thm}{Theorem}
\newtheorem{lem}[thm]{Lemma}
\newtheorem{prop}[thm]{Proposition}
\newtheorem{cor}[thm]{Corollary}
\theoremstyle{definition}
\newtheorem{rem}[thm]{Remark}

\numberwithin{equation}{section}

\def\ie{\emph{i.e.}}
\def\ds{\displaystyle}
\def\:{\colon}
\def\.{\cdot}

\def\<{\left\langle}
\def\>{\right\rangle}
\def\({\left(}
\def\){\right)}
\def\ph#1{\phantom{#1}}
\def\epsilon{\varepsilon}
\def\phi{\varphi}
\def\subset{\subseteq}

\def\leq{\leqslant}
\def\geq{\geqslant}

\def\lra{\longrightarrow}
\def\Lra{\Longrightarrow}
\def\ra{\rightarrow}
\def\bar#1{\overline{#1}}
\def\hat#1{\widehat{#1}}

\def\iso{\cong}
\def\homeq{\simeq}

\DeclareMathOperator{\im}{Im}
\DeclareMathOperator{\Ker}{Ker}\renewcommand{\ker}{\Ker}
\DeclareMathOperator{\ord}{ord}

\def\lie{\mathrm{Lie}}
\def\A{\mathrm{A}}
\def\As{\mathrm{A}^{\mathrm{st}}}
\def\C{\mathbb{C}}

\def\F{\mathbb{F}}

\def\k{\Bbbk}
\def\Q{\mathbb{Q}}

\def\T{\mathcal{T}}
\def\Z{\mathbb{Z}}

\def\einf{E_\infty}

\def\Times_#1{\mathop{\times}_{#1}}
\def\oTimes_#1{\mathop{\otimes}_{#1}}
\def\oPlus_#1{\mathop{\bigoplus}_{#1}}
\def\ideal{\triangleleft}

\DeclareMathOperator{\Cont}{Cont}

\DeclareMathOperator{\Ext}{Ext}

\DeclareMathOperator{\Hom}{Hom}
\DeclareMathOperator{\cHom}{\mathcal{H}om}

\def\CP{\C\mathrm{P}}
\def\id{\mathrm{id}}
\DeclareMathOperator*{\colim}{colim}

\def\HG{\mathrm{H}\Gamma}
\def\cHG{\mathcal{H}\Gamma}

\DeclareMathOperator{\AQ}{AQ}

\def\Ell{E\ell\ell}
\def\emptyset{\varnothing}

\begin{document}
\title[$\Gamma$-cohomology, numerical polynomials, and
$E_\infty$ structures]
{\boldmath$\Gamma$-cohomology of rings of numerical polynomials
and $E_\infty$ structures on $K$-theory}
\author{Andrew Baker}
\address{Mathematics Department,
University of Glasgow, Glasgow G12 8QW, Scotland.}
\email{a.baker@maths.gla.ac.uk}
\urladdr{http://www.maths.gla.ac.uk/$\sim$ajb}
\author{Birgit Richter}
\address{Fachbereich Mathematik der Universit\"at Hamburg,
Bundesstrasse 55, 20146 Hamburg, Germany.}
\email{richter@math.uni-hamburg.de}
\urladdr{http://www.math.uni-hamburg.de/home/richter/}
\thanks{The authors would like to thank the Isaac Newton Institute, 
the Max-Planck-Institut f\"ur Mathematik and the Department of
Mathematics  in Bonn for providing
stimulating environments in which this work was carried out and also
John Greenlees, Alain Jeanneret, Alan Robinson and Stefan Schwede
for helpful comments.}
\begin{abstract}
We investigate $\Gamma$-cohomology of some commutative cooperation algebras
$E_*E$ associated with certain periodic cohomology theories. For $KU$ and
$E(1)$, the Adams summand at a prime $p$, and for $KO$ we show that
$\Gamma$-cohomology vanishes above degree~$1$. As these cohomology groups
are the obstruction groups in the obstruction theory developed by Alan
Robinson we deduce that these spectra admit unique $E_\infty$ structures.
As a consequence we obtain an $E_\infty$ structure for the connective
Adams summand. For the Johnson-Wilson spectrum $E(n)$ with $n\geq1$ we
establish the existence of a unique $E_\infty$ structure for its
$I_n$-adic  completion.
\end{abstract}
\maketitle

\section*{Introduction}

In homotopy theory it is often not sufficient to have homotopy ring
structures on a spectrum in order to construct for instance homotopy
fixed points under a group action or quotient spectra. For this, it
is necessary to have ring structures which are not just given up to
homotopy but where these homotopies fulfil certain coherence conditions.
We will prove the existence and uniqueness of certain $\einf$ structures,
\ie, structures on spectra which encode a coherent homotopy commutative
multiplication.

In~\cite{Ro:Einfty, Ro:unit}, Alan Robinson developed a purely algebraic
obstruction theory for $\einf$ structures on homotopy associative and
commutative ring spectra. The device for deciding whether a spectrum
possesses such a structure is a cohomology theory for commutative algebras,
$\Gamma$-cohomology. When applied to the cooperation algebra $E_*E$ of
a ring spectrum~$E$, the vanishing of these cohomology groups implies the
existence of an $E_\infty$ structure on $E$ which extends
the given homotopy ring structure.

We will apply Robinson's obstruction theory to complex $K$-theory~$KU$,
its $p$-localization $KU_{(p)}$, its Adams summand $E(1)$, with
\[
KU_{(p)}\homeq\bigvee_{i=0}^{p-2}\Sigma^{2i}E(1),
\]
and to real $K$-theory $KO$. Our main topological result is
\begin{Thm}
There are unique $E_\infty$ structures on $KU$, the Adams summand $E(1)$
and $KO$.
\end{Thm}

This process can then be refined to give $E_\infty$ structures on the
connective covers. For the higher Johnson-Wilson spectra $E(n)$ with
$n\geq1$ we will consider the $I_n$-adic completion.

The existence of $\einf$ structures on $KU$ and $KO$ was already known:
in~\cite[Chapter VIII]{May:EinftyRings} $E_\infty$ structures for the
connective versions $ku$ and $ko$ were constructed and the techniques
of~\cite[VIII]{EKMM} lead to $E_\infty$ models for $KU$ and $KO$.
Recently, Joachim~\cite{Joachim} has described such a structure for $KO$
in the context of symmetric spectra. But as far as we know the uniqueness
of these structures has not previously been documented. The existence
and uniqueness for $E(1)$ appears to be new. We also show that the
connective Adams summand $\ell$ at an odd prime~$p$ admits an $E_\infty$
structure; on the $p$-completion $\ell\sphat_p$\; an  $E_\infty$ structure was 
earlier constructed by McClure and Staffeldt~\cite{McC&S}. In
subsequent work we have shown that $E_\infty$ structures on $ku$,
$ko$, $\ell$  and $\ell\sphat_p$ are unique. 

By~\cite{GH,BR&Ro} it is known that the Lubin-Tate spectra $E_n$ have
unique $E_\infty$ structures. In particular, $E_1=KU\sphat_p$\; has a
unique $E_\infty$ structure. The work of Hopkins and Miller~\cite{Rezk}
and Goerss and Hopkins~\cite{GH} establishes $A_\infty$ and $E_\infty$
structures on the Lubin-Tate spectra $E_n$ and identified the homotopy
type of the space of $A_\infty$ (resp.~$E_\infty$) maps between any two
of these spectra. With the help of these results the homotopy action
of the Morava-stabilizer group on $E_n$ could be rigidified to an action
of Morava-stabilizer group on $E_n$ by $A_\infty$ (resp.~$E_\infty$)
maps. 

The existence of unique $E_\infty$ structures on $KU\sphat_p$\; and
$E(1)\sphat_p$\; follow directly from the calculation of continuous
$\Gamma$-cohomology. In~\cite[\S1]{HS}, Hovey and Strickland asked
whether  the $I_n$-adic 
completion of $E(n)$, $\hat{E(n)}$, has $E_\infty$ structures.
\begin{Thm}
For all $n \geq 1$ and all primes $p$, $\hat{E(n)}$ possesses
a unique $E_\infty$ structure.
\end{Thm}

We give an elementary proof of this which relies on Robinson's obstruction
theory~\cite{Ro:Einfty}; the result for $E_n$ then follows using ideas
of~\cite{SVW}. However, so far we have not been able to extend this result
to $E(n)$ itself since the $\Gamma$-cohomology of $E(n)_*E(n)$ appears to
be very non-trivial in positive degrees for $n>1$.

\subsection*{Notation and tools}
All otherwise unspecified tensor products are taken over $\Z$ or a localization
at a prime~$p$, $\Z_{(p)}$. We denote by $\Z_p$, $\Q_p$ and $\hat\Z$
the  rings of
$p$-adic integers, $p$-adic rationals and profinite integers respectively,
while $\hat\Q=\hat\Z\otimes\Q$. For an arbitrary unital commutative ring $R$ we
denote the rationalization of $R$ by $R_\Q$. 

In the body of this paper we will repeatedly use properties of $\Gamma$-cohomology.
For the reader's convenience we recall some of its crucial features here.
\begin{itemize}
\item
Given a commutative $\k$-algebra $A$, $\Gamma$-cohomology of $A$ with 
coefficients in an $A$-module $M$, $\HG^*(A|\k;M)$, is defined 
in~\cite{Ro:Einfty} as the
cohomology of the derived $A$-homomorphisms from a certain chain complex of 
$A$-modules $C^\Gamma(A|\k)_*$ to $M$.
\item
This way of defining $\Gamma$-cohomology ensures that each short exact sequence
of $A$-modules
\[
0\ra M'\lra M\lra M''\ra 0
\]
leads to a long exact sequence of $\Gamma$-cohomology groups. It is also clear
that in good cases there is a universal coefficient spectral sequence
\[
\mathrm{E}_2^{*,*}=\Ext^{*,*}_A(\HG_*(A|\k;A),M)\Lra\HG^*(A|\k;M).
\]
Here $\HG_*(A|\k;A)$ denotes the homology of $C^\Gamma(A|\k)_*$.
\item
In \cite[6.8]{RoWh}, it is shown that $\Gamma$-cohomology vanishes if~$A$
is \'etale over $\k$, and that $\Gamma$-cohomology satisfies 
Flat Base Change and has a Transitivity Sequence. 
\item
Last but not least, if $A$ is a $\k$-algebra and $\Q\subset\k$, then the
$\Gamma$-cohomology of $A$ with coefficients in $M$ coincides with
And\-r\'e\--Qui\-ll\-en cohomology $\AQ^*(A|\k;M)$; for details
see~\cite[6.4]{RoWh}.
\end{itemize}

Convenient and concise sources for the definitions of And\-r\'e\--Qui\-ll\-en
homology and cohomology are~\cite{Loday,Weibel} where they are denoted
$\mathrm{D}_*(\ )$ and $\mathrm{D}^*(\ )$.

The main techniques we use for the $K$-theoretic examples involve the passage
to a continuous version of $\Gamma$-cohomology which we introduce in
Section~\ref{sec:GammCoh}, and a description of the $\Gamma$-cohomology
groups of the relevant cooperation algebras in terms of $\Gamma$-cohomology
groups of colimits of \'etale algebras at each prime (which turn out to vanish) 
and $\Gamma$-cohomology groups of 
algebras which are rationally smooth and can be calculated. At the technical
heart of our arguments for the $K$-theoretic examples is
Theorem~\ref{thm:Z-algebras}.
\section{Linear compactness and cohomology}\label{sec:LinCompact&Cohom}

The results in this section and Section~\ref{sec:GammCoh} will play a
crucial r\^ole in our proof of Corollary~\ref{cor:HG-NumPolys} when we
calculate the $\Gamma$-cohomology of numerical polynomials and more
generally for the proof of Theorem~\ref{thm:Z-algebras}, where
Corollary~\ref{cor:Hom-LinCompact} and Proposition~\ref{prop:Cohom-LinCompact}
will help us to identify certain  $\Gamma$-cohomology groups of complete 
algebras with inverse limits of $\Gamma$-cohomology groups of quotient 
algebras. 

Let $\k$ be a commutative Noetherian ring and let $\mathfrak{m}\ideal\k$
be a maximal ideal. We topologize $\k$ with respect the $\mathfrak{m}$-adic
topology where the open neighbourhoods of $0$ are the ideals
$\mathfrak{m}^k\ideal\k$ for $k\geq0$. In the following we will assume that 
$\k$ is Hausdorff with respect to this $\mathfrak{m}$-adic topology, \ie, that 
\begin{equation}\label{eqn:m-adicHausdorff-k}
\bigcap_{k\geq0}\mathfrak{m}^k=0.
\end{equation}

For each $k\geq0$, $\mathfrak{m}^k$ is a finitely-generated $\k$-module, 
while $\mathfrak{m}^k/\mathfrak{m}^{k+1}$ is a finitely-generated
$\k/\mathfrak{m}$-module which is therefore an Artinian $\k/\mathfrak{m}$-module.

Now let $A$ be a commutative unital $\k$-algebra. The ideals
$\mathfrak{m}^kA=A\mathfrak{m}^k\ideal A$ also generate a topology on~$A$
which is Hausdorff if
\begin{equation}\label{eqn:m-adicHausdorff-A}
\bigcap_{k\geq0}A\mathfrak{m}^k=0.
\end{equation}
Then the unit homomorphism $\k\lra A$ is automatically continuous and if~$A$
is augmented over $\k$ then the augmentation is also continuous. Furthermore,
$(A,\k)$ is a topological algebra over the topological ring $\k$. We say that
$(A,\k)$ is Hausdorff if both~\eqref{eqn:m-adicHausdorff-k}
and~\eqref{eqn:m-adicHausdorff-A} hold.

If $(A,\k)$ is a Hausdorff topological algebra, then the \emph{$\mathfrak{m}$-adic
completion} of $(A,\k)$ is $(A\sphat_{\mathfrak{m}},\k\sphat_{\mathfrak{m}})$,
where
\[
A\sphat_{\mathfrak{m}}=\lim_kA/\mathfrak{m}^k\!A,
\quad
\k\sphat_{\mathfrak{m}}=\lim_k\k/\mathfrak{m}^k.
\]
We say that $(A,\k)$ is \emph{$\mathfrak{m}$-adically complete} if
$A\sphat_{\mathfrak{m}}=A$ and $\k\sphat_{\mathfrak{m}}=\k$. When
$\mathfrak{m}$ is clear from the context we will sometimes simplify
notation by writing
$(\hat{A},\hat{\k})=(A\sphat_{\mathfrak{m}},\k\sphat_{\mathfrak{m}})$.
If $(A,\k)$ is augmented over $\k$ then
$(A\sphat_{\mathfrak{m}},\k\sphat_{\mathfrak{m}})$ is augmented over
$\k\sphat_{\mathfrak{m}}$.

Recall from~\cite{Bourbaki:CommAlg,Jensen:LNM,Zelinsky} the property
of a topological left module $M$ over a topological algebra $(A,\k)$
being \emph{linearly compact}:
\begin{itemize}
\item
If $\{x_\lambda+N_\lambda\}_{\lambda\in\Lambda}$ is a collection of cosets
of closed submodules $N_\lambda\subset M$ such that every finite collection
of the cosets $x_\lambda+N_\lambda$ has non-trivial intersection, then
\[
\bigcap_{\lambda\in\Lambda}x_\lambda+N_\lambda\neq\emptyset.
\]
\end{itemize}
We will make repeated use of the following vanishing result
of~\cite[p.~57, th\'eor\`eme~7.1]{Jensen:LNM} for the higher derived
functors of the inverse limit for inverse systems of linearly compact
modules.
\begin{thm}\label{thm:LInComp->Rlim=0}
Let $\{M_i\}_{i \in I}$ be an inverse system of linearly compact $A$-modules
and continuous\/ $A$-linear maps. Then for all $s>0$ we have
\[
{\lim_i}^sM_i=0.
\]
\end{thm}

Recall that a topological $\k$-module $M$ is \emph{topologically free on
a countable basis $\{b_i\}_{i\geq1}$} if for each element $x\in M$ and
$k\geq1$, in $M/\mathfrak{m}^kM$ considered as an $\k/\mathfrak{m}^k$-module,
there is a unique (finite) expansion
\[
\bar{x}=\sum_{i\geq1}\bar{r}_i \bar{b}_i
\]
with $\bar{r}_i\in\k/\mathfrak{m}^k$ and where $\bar{b}_i\in M/\mathfrak{m}^kM$
is the residue class of $b_i$. As a consequence, $x$ has a unique expansion
as a limit sum
\[
x=\sum_{i\geq1}t_ib_i
\]
where $t_i\ra0$ as $i\ra\infty$; this means that for each $k$, there is
an $n_k$ such that for $i>n_k$ we have $t_i\in\mathfrak{m}^k$. The linear
topology on $M$ has basic open neighbourhoods of~$0$ of the form
$\mathfrak{m}^kM$. Now the Noetherian condition on~$\k$ implies that
\begin{equation}\label{eqn:LinTop-m^kM}
\mathfrak{m}^kM=
\{\sum_{i\geq1}t_ib_i:\text{$t_i\in\mathfrak{m}^k$, $t_i\ra0$ as $i\ra\infty$}\}.
\end{equation}

For two topological left $R$-modules $L$ and $M$ over a commutative topological
ring $R$, we let
\[
\cHom_{R}(L,M)\subset\Hom_{R}(L,M)
\]
be the submodule of continuous $R$-module homomorphisms.

\begin{prop}\label{prop:Hom-LinCompact}
Suppose that $M$ is a finitely generated $\k$-module which is
complete and Hausdorff with respect to the $\mathfrak{m}$-adic topology.
If~$L$ is a $\k$-module which is complete with respect to the
$\mathfrak{m}$-adic topology and topologically free on a countable
basis then $\cHom_\k(L,M)$ is a linearly compact $\k$-module.
\end{prop}
\begin{proof}
Assume first that $L$ is Hausdorff with respect to the $\mathfrak{m}$-adic
topology. 
Note that
\[
\cHom_\k(L,M)=\cHom_\k(L,\lim_kM/\mathfrak{m}^kM)
=\lim_k\cHom_\k(L,M/\mathfrak{m}^kM).
\]
If $\{b_j\}_{j\geq1}$ is a topological basis for $L$, then using the
Noetherian condition on $\k$ we find that the basic neighbourhoods
of~$0$ in $L$ are the submodules $\mathfrak{m}^kL\subset L$. From
this we find that
\begin{align*}
\cHom_\k(L,M/\mathfrak{m}^kM)
&=\cHom_{\k/\mathfrak{m}^k}(L/\mathfrak{m}^kL,M/\mathfrak{m}^kM) \\
&=\Hom_{\k/\mathfrak{m}^k}(L/\mathfrak{m}^kL,M/\mathfrak{m}^kM) \\
&=\prod_{j\geq1}\Hom_{\k/\mathfrak{m}^k}((\k/\mathfrak{m}^k)b_j,M/\mathfrak{m}^kM).
\end{align*}
But
\[
\Hom_{\k/\mathfrak{m}^k}((\k/\mathfrak{m}^k)b_j,M/\mathfrak{m}^kM)
=M/\mathfrak{m}^kM
\]
and this is Artinian, hence linearly compact. This in turn implies that
the final product above is also linearly compact. The claim now follows
since $\lim_k\cHom_\k(L,M/\mathfrak{m}^kM)$ is a closed subspace of the
product
\[
\prod_{j\geq1}\Hom_{\k/\mathfrak{m}^k}((\k/\mathfrak{m}^k)b_j,M/\mathfrak{m}^kM).
\]

Now we consider what happens when $L$ is not necessarily Hausdorff.
In this case, Nakayama's Lemma implies that for any $f\in\cHom_\k(L,M)$
we have
\[
f \Big(\bigcap_{k\geq1}\mathfrak{m}^kL\Big)=0.
\]
Hence such an $f$ factors through the quotient
\[
L_0=L\Big/\!\bigcap_{k\geq1}\mathfrak{m}^kL,
\]
so we might as well replace $L$ by this Hausdorff quotient. Then
we have
\begin{equation}\label{eqn:Hom-LinCompact-Hausdorff}
\cHom_\k(L,M)=\cHom_\k(L_0,M).
\qedhere
\end{equation}
\end{proof}

We can apply this to prove the following.
\begin{cor}\label{cor:Hom-LinCompact}
Suppose further that $A$ is a topological $\k$-algebra with respect
to the $\mathfrak{m}$-adic topology inherited from $\k$ and that $L$
and $M$ are topological $A$-modules. Then
$\cHom_A(L,M)\subset\cHom_\k(L,M)$ is a closed $\k$-submodule.
Hence $\cHom_A(L,M)$ is linearly compact.
\end{cor}
\begin{proof}
Again we first consider the case where $L$ is Hausdorff. 
The two continuous action maps
\[
A\otimes_\k\cHom_\k(L,M)\lra\cHom_\k(L,M)
\]
given by
\[
a\otimes f\longmapsto af,\quad a\otimes f\longmapsto f(a(-))
\]
are equalized on $\cHom_A(L,M)$, so this is a closed subset
of $\cHom_\k(L,M)$.
For $L$ not Hausdorff we see as above that $\cHom_A(L,M) \cong 
\cHom_A(L_0,M)$. 
\end{proof}
Note that if $A$ is not necessarily Hausdorff, then setting
\[
A_0=A\Big/\bigcap_{k\geq1}\mathfrak{m}^kA
\]
we have
\begin{equation*} 
\cHom_A(L,M)=\cHom_{A_0}(L_0,M).
\end{equation*}
\begin{prop}\label{prop:Cohom-LinCompact}
Let $(C^*,\delta)$ be a cochain complex of linearly compact and Hausdorff\/
$\k$-modules where for each $n$, the coboundary $\delta^n\:C^n\lra C^{n+1}$
is continuous. Then for each $n$, $\mathrm{H}^n(C^*,\delta)$ is linearly
compact.
\end{prop}
\begin{proof}
Since each $C^n$ is linearly compact and Hausdorff, the submodules
$\im\delta^{n-1}$ and $\ker\delta^n$ of $C^n$ are both closed. Therefore
\[
\mathrm{H}^n(C^*,\delta)=\ker\delta^n/\im\delta^{n-1}
\]
is also linearly compact.
\end{proof}

\section{Continuous $\Gamma$-cohomology}\label{sec:GammCoh}

In this section we discuss some technical issues related to our calculations
of $\Gamma$-co\-hom\-ol\-ogy later in the paper. Continuous cohomology of 
profinite groups is described in~\cite{Shatz,Weibel}; for analogues appearing 
in topology see~\cite{AB:Ainfty,AB:HopfAlgSS}; our present theory is modelled 
closely on the presentations in those references. 

Let $\k$ be a commutative Noetherian ring and let $\mathfrak{m}\ideal\k$
be a maximal ideal. We topologize $\k$ with the $\mathfrak{m}$-adic
topology. Let $M$ be a topological left module over a topological algebra 
$(A,\k)$. In practise, we will usually consider the $\mathfrak{m}$-adic 
topology on $M$.

In the following we consider $\Gamma$-cohomology of $A$ with coefficients 
in  $M$, $\HG^*(A|\k;M)$. In~\cite{Ro:Einfty,RoWh} $\Gamma$-cohomology 
is defined using a cochain complex $\Hom_A(C^\Gamma(A|\k)_*,M)$. Here 
$C^\Gamma(A|\k)_*$ is the chain complex whose degree $n$-part is the 
left $A$-module 
$$
C^\Gamma(A|\k)_n=\bigoplus_{s+r=n} \lie(s+1)^* \otimes 
\k[\Sigma_{s+1}]^{\otimes r} \otimes_\k A^{\otimes (s+2)}, 
$$
where $\lie(s+1)^*$ is the $\k$-linear dual of the $(s+1)$-st term of the
Lie operad and $\Sigma_\ell$ denotes the symmetric group on $\ell$ letters. 
In particular, $\lie(s+1)^*$ and $\k[\Sigma_{s+1}]$ are finite-rank free 
$\k$-modules. Topologising $C^\Gamma(A|\k)_*$ with the $\mathfrak{m}$-adic 
topology, we can introduce the subcomplex
\[
\cHom_{A}(C^\Gamma(A|\k)_*,M)\subset\Hom_{A}(C^\Gamma(A|\k)_*,M)
\]
of continuous cochains whose cohomology $\cHG^*(A|\k;M)$ we call the
\emph{continuous $\Gamma$-cohomology} of $A$ with coefficients in~$M$.
Note that continuous cochains can be expressed as the inverse limit
$$
\cHom_{A}(C^\Gamma(A|\k)_*,M) = 
\lim_k\Hom_{A/\mathfrak{m}^k\!A} 
(C^\Gamma(A/\mathfrak{m}^k\!A|\k/\mathfrak{m}^k)_*, M/\mathfrak{m}^kM) 
$$
if $M$ is Hausdorff. 

\begin{lem} \label{lem:2cake}
If $(A, \k)$ is a topological algebra as above whose completion $\hat{A}$ 
is countably free on a topological basis as a $\hat{\k}$-module and if $M$
is an $A$-module such that $M/\mathfrak{m}^kM$ is Artinian and $M$ is complete,
then $\Hom_{A}(C^\Gamma(A|\k)_n,M)$ is linearly compact in each degree $n$. 
\end{lem}
\begin{proof}
As $M$ is complete, 
\begin{align} \label{eq:milnor}
\notag \Hom_{A}(C^\Gamma(A|\k)_n,M) & = \Hom_{A}(C^\Gamma(A|\k)_n, 
\lim_k M/\mathfrak{m}^kM) \\
 & = \lim_k \Hom_{A}(C^\Gamma(A|\k)_n, M/\mathfrak{m}^kM). 
\end{align}
For any fixed $k$, the homomorphisms from $C^\Gamma(A|\k)_n$ to the quotient 
$M/\mathfrak{m}^kM$ factor through  
\begin{align*}
C^\Gamma(A/\mathfrak{m}^kA|\k/\mathfrak{m}^k)_n & = 
\bigoplus_{s+r=n} \lie(s+1)^* \otimes \k/\mathfrak{m}^k[\Sigma_{s+1}]^{\otimes r} 
\otimes (A/\mathfrak{m}^kA)^{\otimes (s+2)} \\
& = C^\Gamma(\hat{A}/\mathfrak{m}^k\hat{A}|\hat{\k}/\mathfrak{m}^k)_n.
\end{align*}
As we assumed that $\hat{A}$ is countably free on a topological basis, the 
quotient $\hat{A}/\mathfrak{m}^k\hat{A}$ is free on a countable basis.
Therefore 
for any homomorphism $f$ in $\Hom_{A}(C^\Gamma(A|\k)_n,M)$ there is an $j$ 
such that $f$ factors over the finitely-generated free submodule 
$C(j)_n \subset C^\Gamma(\hat{A}/\mathfrak{m}^k\hat{A}|\hat{\k}/\mathfrak{m}^k)_n$ 
spanned by the first $j$ generators of 
$C^\Gamma(\hat{A}/\mathfrak{m}^k\hat{A}|\hat{\k}/\mathfrak{m}^k)_n$. Therefore 
\begin{align*}
\Hom_{A}(C^\Gamma(A|\k)_n,M) & = \lim_k \Hom_{A}(C^\Gamma(A|\k)_n, M/\mathfrak{m}^kM) \\
 & = \lim_k \Hom_{\hat{A}}(C^\Gamma(\hat{A}/\mathfrak{m}^k\hat{A}|\hat{\k}/\mathfrak{m}^k)_n, 
M/\mathfrak{m}^kM) \\
 & = \lim_k \lim_j \Hom_{\hat{A}}(C(j)_n, M/\mathfrak{m}^kM) \\
 & = \lim_k \lim_j \prod_{i=1}^j M/\mathfrak{m}^kM.
\end{align*}
Therefore $\Hom_{A}(C^\Gamma(A|\k)_n,M)$ is a limit of Artinian modules and 
thus linearly compact. 
\end{proof}
Notice that the above inclusion of complexes induces a forgetful homomorphism
\begin{equation}\label{eqn:ContsCohom-Cohom}
\rho\:\cHG^*(A|\k;M)\lra\HG^*(A|\k;M).
\end{equation}
Recall that if  $M$ is $\mathfrak{m}$-adically complete and Hausdorff, 
there is a short exact sequence
\begin{equation} \label{eqn:sesforM}
0\ra M\lra\prod_kM/\mathfrak{m}^kM\xrightarrow{\id-\sigma}
\prod_kM/\mathfrak{m}^kM\ra0,
\end{equation}
where $\sigma$ is the shift-reduction map. From (\ref{eq:milnor}) and 
(\ref{eqn:sesforM}) we deduce a Milnor exact 
sequence relating 
$\cHG^*$ to ordinary $\Gamma$-cohomology for complete  coefficient modules. 
For a similar result see~\cite{AB:Ainfty}.
\begin{prop}\label{prop:ContGH-fgmodule}
As above, let $\k$ be Noetherian with maximal ideal $\mathfrak{m}$. Let $M$ 
be a complete Hausdorff topological module over $\hat{A}$ which is
finitely-generated over $\hat{\k}$. Then for each $n$ there is a short exact 
sequence
\begin{multline*}
0\ra{\lim_k}^1\HG^{n-1}(A/\mathfrak{m}^k\!A|\k/\mathfrak{m}^k;M/\mathfrak{m}^kM)
\lra\cHG^n(\hat{A}|\hat{\k};M) \\
\lra\lim_k\HG^n(A/\mathfrak{m}^k\!A|\k/\mathfrak{m}^k;M/\mathfrak{m}^kM)\ra0.
\end{multline*}
\end{prop}

This leads to some useful calculational results, versions of which have previously
appeared in~\cite{Rezk,BR&Ro}. Notice that for any $\hat{A}$-module~$M$ and $k\geq1$,
there is a natural reduction homomorphism
\begin{equation}\label{eqn:Naturality}
\HG^n(\hat{A}|\hat{\k};M)\lra\HG^n(A|\k;M)
\lra\HG^n(A/\mathfrak{m}^kA|\k/\mathfrak{m}^k;M/\mathfrak{m}^kM),
\end{equation}
compatible with respect to different values of $k$. 
In turn there is a homomorphism
\begin{equation}\label{eqn:Naturality-lim}
\HG^n(\hat{A}|\hat{\k};M)\lra\HG^n(A|\k;M)
\lra
\lim_k\HG^n(A/\mathfrak{m}^k\!A|\k/\mathfrak{m}^k;M/\mathfrak{m}^kM).
\end{equation}

The following result was inspired by~\cite[lemma~15.6]{Rezk}. 
\begin{cor}\label{prop:ContGH-MilnorSeq}
Let $M$ be an $\hat{A}$-module which is complete and Hausdorff with respect 
to the $\mathfrak{m}$-adic topology and finitely-generated as a $\hat{\k}$-module. 
Let $\hat{A}$ be countably free on a topological basis. Then the natural homomorphism 
$\rho$ induces an isomorphism
\[
\cHG^*(\hat{A}|\hat{\k};M)\iso\HG^*(\hat{A}|\hat{\k};M). 
\]
In addition 
\[
{\lim_k}^1\HG^{n-1}(A/\mathfrak{m}^k\!A|\k/\mathfrak{m}^k; 
M/\mathfrak{m}^kM)=0
\]
and the natural homomorphisms induce isomorphisms
\[
\HG^*(A|\k;M)\iso\HG^*(\hat{A}|\hat{\k};M) \iso 
\cHG^*(\hat{A}|\hat{\k};M)
\iso\lim_k\HG^n(A/\mathfrak{m}^k\!A|\k/\mathfrak{m}^k;M/\mathfrak{m}^kM).
\]
\end{cor}
\begin{proof}
Using the naturality provided by~\eqref{eqn:Naturality} we obtain a diagram 
of short exact sequences from the Milnor exact sequence of 
Proposition~\ref{prop:ContGH-fgmodule} into the one for $\HG^n(\hat{A}|\hat{\k};M)$. 
As the homomorphisms at either end are identities, the natural map
$\cHG^*(\hat{A}|\hat{\k};M)\lra\HG^*(\hat{A}|\hat{\k};M)$ is an isomorphism.

Under the assumptions, the cochain complex for $\Gamma$-cohomology is linearly 
compact in each degree, by Lemma~\ref{lem:2cake}. Hence by 
Proposition~\ref{prop:Cohom-LinCompact}, $\Gamma$-cohomology is also linearly 
compact in each degree. Therefore 
\[
{\lim_k}^1\HG^{n-1}(A/\mathfrak{m}^k|\k/\mathfrak{m}^k;M/\mathfrak{m}^kM)=0.
\]
So in this case
\[
\cHG^*(\hat{A}|\hat{\k};M)\iso\HG^*(\hat{A}|\hat{\k};M)
\iso\lim_k\HG^n(A/\mathfrak{m}^k\!A|\k/\mathfrak{m}^k;M/\mathfrak{m}^kM).
\qedhere 
\]
\end{proof}
\begin{rem}\label{rem:cHH}
Analogous ideas apply to Hochschild cohomology for which a continuous
version appears in~\cite{AB:Ainfty}.
\end{rem}

The following result which will be used in Section~\ref{sec:HG-NumPolys}.
\begin{prop}\label{prop:LESk-k^} 
Let $\k$ be Noetherian with maximal ideal $\mathfrak{m}$ and let $\hat{\k} = 
\lim_k \k/\mathfrak{m}^k$. For any $\k$-algebra $A$ with the $\mathfrak{m}$-adic 
topology for which $\hat{A}$ is topologically free on a countable basis there 
is a long exact sequence
\begin{multline*}
\cdots\lra\HG^{n-1}(A|\k;\hat{\k}/\k)\lra
\HG^n(A|\k;\k)\lra\HG^n(A|\k;\hat{\k})\lra
\HG^n(A|\k;\hat{\k}/\k)\\
\lra\HG^{n+1}(A|\k;\k)\lra\cdots.
\end{multline*}
Making use of the isomorphism $\HG^n(\hat{A}|\hat{\k};\hat{\k}) \cong 
\HG^n(A|\k;\hat{\k})$, we also obtain an analogous exact sequence for 
$\HG^n(\hat{A}|\hat{\k};\hat{\k})$. 
\end{prop}
\begin{proof}
The short exact sequence
\[
0\ra\k\lra\hat{\k}\lra\hat{\k}/\k\ra0
\]
of coefficients together with the last isomorphism from 
Proposition~\ref{prop:ContGH-MilnorSeq} yields this long exact sequence.
\end{proof}

Finally, we record a result on the $\Gamma$-(co)homology of formally
\'etale algebras that we will repeatedly use. We call an algebra
\emph{formally \'etale} if it is a colimit of \'etale algebras.
\begin{lem}\label{lem:HG-Formallyetale}
If $(A,\k)$ is a formally \'etale algebra then for any $A$-module $M$,
\[
\HG_*(A|\k;M)=0=\HG^*(A|\k;M).
\]
\end{lem}
\begin{proof}
By~\cite[theorem 6.8~(3)]{RoWh}, $\Gamma$-homology and cohomology vanishes
for \'etale algebras. Also, $\Gamma$-homology commutes with colimits. Hence
if $A=\ds\colim_{\ph{lc}r}A_r$ with $A_r$ \'etale, for any $A$-module $M$
we have
\[
\HG_*(A|\k;M)=\colim_{\ph{lc}r}\HG_*(A_r|\k;M)=0.
\]
The universal coefficient spectral sequence
\[
\mathrm{E}_2^{*,*}=\Ext^{*,*}_A(\HG_*(A|\k;A),M)\Lra\HG^*(A|\k;M),
\]
has trivial $\mathrm{E}_2$-term, therefore $\HG^*(A|\k;M)=0$.
\end{proof}

\section{Rings of numerical polynomials}\label{sec:NumPolys}

We need to describe some properties of rings of numerical polynomials
which appeared in a topological setting in~\cite{AHS,BCRS} and we follow
these sources in our discussion. As topological motivation, we remark that
$\A$ can be identified with $KU_0\CP^\infty$ and $\As$ with $KU_0KU$ and
we will calculate the $\Gamma$-cohomology of $KU_0KU$ later. By definition,
\begin{align*}
\A&=\{f(w)\in\Q[w]:\forall n\in\Z,\;f(n)\in\Z\}, \\
\As&=\{f(w)\in\Q[w,w^{-1}]:\forall n\in\Z-\{0\},\;f(n)\in\Z[1/n]\}.
\end{align*}
are the rings of \emph{numerical} and \emph{stably numerical} polynomials
(over $\Z$). If $x,y$ are indeterminates, we can work in any of the rings
$\A[x,y]$, $\As[x,y]$ or $\Q[w,w^{-1}][x,y]$.

We will make use of the binomial coefficient functions
\[
c_n(w)=\binom{w}{n}=\frac{w(w-1)\cdots(w-n+1)}{n!}\in A\subset\Q[w]
\]
which can be encoded in the generating function
\[
(1+x)^w=\sum_{n\geq0}c_n(w)x^n\in\A[x]\subset\Q[w][x].
\]
Notice that this satisfies the formal identity
\begin{equation}\label{eqn:BinomIdentity-GenFn}
(1+x)^w(1+y)^w=(1+(x+y+xy))^w.
\end{equation}
Thus we have
\begin{equation}\label{eqn:BinomIdentity-Explicit}
c_m(w)c_n(w)=\binom{m+n}{m}c_{m+n}(w)+(\text{\rm terms of lower degree})
\quad \text{ for }  m,n\geq0.
\end{equation}
\begin{thm}\label{thm:NumPolys-Facts}
\emph{\cite{Ad-Cl, BCRS}}
\begin{enumerate}
\item[(a)]
$\A$ is a free $\Z$-module with a basis consisting of the $c_n(w)$
for $n\geq0$.
\item[(b)]
$\As$ is the localization $\As=\A[w^{-1}]$ and it is a free $\Z$-module
on a countable basis.
\end{enumerate}
\end{thm}

Describing explicit $\Z$-bases for $\As$ is a non-trivial task, 
see~\cite{CCW,Johnson:K*K-Basis}. On the other hand, the multiplicative
structure of the $\Z$-algebra $\As$ is in some ways more understandable. Our
next result describes some generators for $\As$.
\begin{thm}\label{thm:NumPolys-Mult}
\emph{\cite{AHS,BCRS}}
\begin{enumerate}
\item[(a)]
The $\Z$-algebra $\A$ is generated by the elements $c_m(w)$ with $m\geq1$
subject to the relations of~\eqref{eqn:BinomIdentity-Explicit}.
\item[(b)]
The $\Z$-algebra, $\As$ is generated by the elements $w^{-1}$ and $c_m(w)$
with $m\geq1$.
\item[(c)]
We have
\[
\A\otimes\Q=\Q[w],\quad\As\otimes\Q=\Q[w,w^{-1}].
\]
\end{enumerate}
\end{thm}

For the localizations of the rings $\A$ and $\As$ at any prime $p$ 
we have
\begin{subequations}\label{eqn:Local}
\begin{align}
\A_{(p)}&=\{f(w)\in\Q[w]:\forall u\in\Z_{(p)},\;f(u)\in\Z_{(p)}\},
\label{eqn:A-Local}\\
\As_{(p)}&=\{f(w)\in\Q[w]:\forall u\in\Z_{(p)}^\times,\;f(u)\in\Z_{(p)}\}.
\label{eqn:As-Local}
\end{align}
\end{subequations}
\begin{thm}\label{thm:NumPolys-Local}
{\rm(}\cite{Ad-Cl}, \cite[prop.~2.5]{AB:p-adic}, \cite{BCRS}{\rm)}
\begin{enumerate}
\item[(a)]
$\A_{(p)}$ is a free $\Z_{(p)}$-module with a basis consisting of
the monomials in the binomial coefficient functions
\[
w^{r_0}c_p(w)^{r_1}c_{p^2}(w)^{r_2}\cdots c_{p^\ell}(w)^{r_\ell},
\]
where $r_k=0,1,\ldots,p-1$.
\item[(b)]
The $\Z_{(p)}$-algebra $\A_{(p)}$ is generated by the elements $c_{p^m}(w)$
with $m\geq0$ subject to relations of the form
\[
c_{p^m}(w)^p-c_{p^m}(w)=p\,d_{m+1}(w),
\]
where $d_{m+1}(w)\in\A_{(p)}$ has $\deg d_{m+1}(w)=p^{m+1}$. In fact the
monomials
\[
w^{r_0}d_1(w)^{r_1}d_2(w)^{r_2}\cdots d_\ell(w)^{r_\ell},
\]
where $r_k=0,1,\ldots,p-1$, form a basis of $\A_{(p)}$ over $\Z_{(p)}$ and
are subject to multiplicative relations of the form
\[
d_m(w)^p-d_m(w)=p\,d'_{m+1}(w),
\]
where $\deg d'_{m+1}(w)=p^{m+1}$.
\item[(c)]
$\As_{(p)}$ is the localization $\As_{(p)}=\A_{(p)}[w^{-1}]$ and it is a
free $\Z_{(p)}$-module on a countable basis.
\item[(d)]
The $\Z_{(p)}$-algebra, $\As_{(p)}$ is generated by the elements $w$
and $e_m(w)\in\As_{(p)}$ for $m\geq1$ defined recursively by
\[
w^{p-1}-1=pe_1(w),\qquad e_m(w)^p-e_m(w)=pe_{m+1}(w)\quad \text{ for }
m\geq1.
\]
\end{enumerate}
\end{thm}
\begin{cor}\label{cor:NumPolys-Local}
Let $p$ be a prime. 
\begin{enumerate}
\item[(a)]
As $\F_p$-algebras,
\begin{align*}
\A/p\A&=\F_p[c_{p^m}(w):m\geq0]/(c_{p^m}(w)^p-c_{p^m}(w):m\geq0), \\
\As/p\As&=\F_p[w,e_m(w):m\geq0]/(w^{p-1}-1,e_m(w)^p-e_m(w):m\geq1).
\end{align*}
Hence these algebras are formally \'etale over $\F_p$.
\item[(b)]
For $n\geq1$, $\A/p^n\A$ and $\As/p^n\As$ are formally \'etale
over $\Z/p^n$.
\item[(c)]
The $p$-adic completions $\A_p=\lim_{n}\A/p^n\A$ and
$\As_p=\lim_{n}\As/p^n\As$ are formally \'etale over $\Z_p$.
\item[(d)]
$\A_p$ and $\As_p$ are free topological $\Z_p$-modules on countable bases.
Therefore they are both compact and Hausdorff.   
\end{enumerate}
\end{cor}
\begin{proof}
Parts (b) and hence (c) can be proved by induction on $n\geq1$ using
the \emph{infinite-dimensional Hensel lemma} of~\cite[3.9]{AB:Ainfty}.
The case $n=1$ is immediate from~(a). Suppose that we have found a
sequence of elements $s_0,s_1,\ldots,s_k,\ldots\in\A_{(p)}$ satisfying
\[
s_m^p-s_m\equiv0\pmod{p^n}\quad \text{ for } m\geq0.
\]
Taking $s'_m=s_m+(s_m^p-s_m)$ we find that
\begin{align*}
{s'_m}^p-s'_m&=(s_m+(s_m^p-s_m))^p-(s_m+(s_m^p-s_m)) \\
&\equiv s_m^p-(s_m+(s_m^p-s_m))\pmod{p^{n+1}} \\
&=0.
\end{align*}
Hence for every $n$ we can inductively produce such elements
$s_{n,m}\in\A_{(p)}$ for which
\begin{align*}
\A/p^n\A&=\Z/p^n[s_{n,m}:m\geq0]/(s_{n,m}^p-s_{n,m}:m\geq0) \\
&=\bigotimes_{m\geq0}\Z/p^n[s_{n,m}]/(s_{n,m}^p-s_{n,m}).
\end{align*}
Now passing to $p$-adic limits we obtain elements
\[
s_m=\lim_{n\ra\infty}s_{n,m}\in\A_p
\]
for which
\[
s_m^p-s_m=0.
\]
In these cases we obtain for the module of K\"ahler differentials
\[
\Omega^1_{(\A/p^n\A)/\Z/p^n}=0=\Omega^1_{\A_p/\Z_p}.
\]
Part (d) is related to Mahler's Theorem and a suitable exposition of
this can be found in \cite{AB:p-adic}.
\end{proof}

There are two natural choices of augmentation for $\A$, namely
evaluation at $0$ or $1$,
\begin{align*}
\epsilon_+\:\A\lra\Z;\quad&\epsilon_+f(w)=f(0), \\
\epsilon_\times\:\A\lra\Z;\quad&\epsilon_\times f(w)=f(1).
\end{align*}
For our purposes, the latter augmentation will be used. Notice
that there is a ring automorphism
\begin{equation}\label{eqn:A-Auto}
\phi\:\A\lra\A;\quad\phi f(w)=f(w+1)
\end{equation}
for which $\epsilon_+\phi=\epsilon_\times$, so these augmentations
are not too dissimilar.

\section{The ring of $\Z/(p-1)$-invariants in $\As_{(p)}$}
\label{sec:Z/(p-1)invariants}

In this section, $p$ always denotes an \emph{odd} prime. The case
of $p=2$ is related to $KO$ and the work of Section~\ref{sec:KO}.

Since polynomial functions $\Z_{(p)}^\times\lra\Q$ are continuous
with respect to the $p$-adic topology they extend to continuous
functions $\Z_p^\times\lra\Q_p$; such functions which also map
$\Z_{(p)}^\times$ into $\Z_{(p)}$ give continuous functions
$\Z_p^\times\lra\Z_p$. Hence we can regard $\As_{(p)}$ as a subring
of $\Q_p[w,w^{-1}]$ which in turn can be viewed as a space of
continuous functions on the $p$-adic units $\Z_p^\times$. For $p \geq
3$ there  is a splitting of topological groups
\[
\Z_p^\times\iso\Z/(p-1)\times(1+p\Z_p),
\]
where $\Z/(p-1)$ identifies with a subgroup generated by a primitive
$(p-1)$-st root of unity $\zeta$. There is also a bicontinuous
isomorphism $1+p\Z_p\iso\Z_p$. 

For an odd prime $p$, the group $\<\zeta\>\iso\Z/(p-1)$ acts continuously
on $\Q_p[w,w^{-1}]$ by
\[
\zeta\.f(w)=f(\zeta w)
\]
and it is immediate that this action sends elements of $\As_{(p)}$ to
continuous functions $\Z_p^\times\lra\Z_p$. It then makes sense to ask
for the subring of $\As_{(p)}$ fixed by this action, ${}^\zeta\As_{(p)}$.
We will relate this subring to the algebra of cooperations of the Adams
summand in Proposition \ref{prop:KU*KU}.

Recall the elements $e_m(w)$ of Theorem~\ref{thm:NumPolys-Local}(d).
We will write $\bar{e}_m(w)$ for $w^{-1}e_m(w)$.
\begin{prop}\label{prop:Intsubring}
As a $\Z_{(p)}$-algebra, ${}^\zeta\As_{(p)}$ is generated by the
elements $w^{p-1}$ and $\bar{e}_m(w)$ for $m\geq1$.
\end{prop}
\begin{proof}
It is clear that
\[
{}^\zeta\Q[w,w^{-1}]=\Q[w^{p-1},w^{-(p-1)}].
\]
Also, by construction of the $e_m(w)$,
\[
\bar{e}_m(w)\in{}^\zeta\As_{(p)}\subset\Q[w^{p-1},w^{-(p-1)}].
\]
Consider the multiplicative idempotent
\[
E_\zeta\:\Q[w,w^{-1}]\lra\Q[w,w^{-1}];
\quad
E_\zeta f(w)=\frac{1}{p-1}\sum_{r=1}^{p-1}f(\zeta^rw).
\]
Then we have
\[
{}^\zeta\As_{(p)}=E_\zeta\As_{(p)}.
\]
Each element $f(w)\in\Q[w,w^{-1}]$ has the form
\[
f(w)=f_0(w^{p-1})+wf_1(w^{p-1})+\cdots+w^{p-2}f_{p-2}(w^{p-1}), 
\]
where $f_k(x)\in\Q[x,x^{-1}]$, hence
\[
E_\zeta f(w)=f_0(w^{p-1}).
\]
From this it follows that ${}^\zeta\As_{(p)}$ is generated
as an $\Z_{(p)}$-algebra by the stated elements.
\end{proof}
\begin{cor}\label{cor:Intsubring}
The following hold.
\begin{enumerate}
\item[(a)]
As $\F_p$-algebras,
\[
{}^\zeta\As_{(p)}/p({}^\zeta\As_{(p)})=
\F_p[w,\bar{e}_m(w):m\geq1]/(w^{p-1}-1,\bar{e}_m(w)^p-\bar{e}_m(w):m\geq1).
\]
Hence this algebra is formally \'etale over $\F_p$.
\item[(b)]
For $n\geq1$, ${}^\zeta\As_{(p)}/p^n({}^\zeta\As_{(p)})$ is formally
\'etale over $\Z/p^n$.
\item[(c)]
The $p$-adic completion
${}^\zeta\As_p=\lim_{n}{}^\zeta\As_{(p)}/p^n({}^\zeta\As)_{(p)}$
is formally \'etale over $\Z_p$.
\end{enumerate}
\end{cor}

\section{The $\Gamma$-cohomology of numerical polynomials}
\label{sec:HG-NumPolys}

Recall that $\hat{\Z}/\Z$ and $\Z_p/\Z_{(p)}$ for any prime~$p$ are
torsion-free divisible groups, so they are both $\Q$-vector spaces
which have the same cardinality and (uncountable) dimensions; thus
they are isomorphic. Similarly, we have
$\hat{\Z}/\Z\iso\hat{\Q}/\Q$ and $\Z_p/\Z_{(p)}\iso\Q_p/\Q$.

From now on, we will always use the augmentations $\epsilon_\times\:\A\lra\Z$
and $\epsilon_\times\:\As\lra\Z$ and their analogues for the $p$-localized
versions. In calculating the $\Gamma$-cohomology of $\A$, we would obtain
the same result using $\epsilon_+$ because of the existence of the
automorphism~$\phi$ of~\eqref{eqn:A-Auto}.
\begin{thm}\label{thm:Z-algebras}
Let $R$ be an augmented commutative $\Z$-algebra. 
Assume that, at each prime $p$, the $p$-completion $R_p$ is topologically free 
on a countable basis. 
Suppose that for all primes $p$ and  $k\geq 1$, $R/p^k$ is a 
formally \'etale algebra over $\Z/p^k$. 
Then for all $s \geq 0$,
\[
\HG^s(R|\Z;\Z) \cong \HG^{s-1}(R_\Q|\Q;\Q).
\]
\end{thm}
\begin{proof}
For each natural number~$n$, we may write 
\[
n=\prod_{p}p^{\ord_p n}
\]
where the product is taken over all primes~$p$. The Chinese Remainder
Theorem provides splittings
\begin{subequations}\label{eqn:ProfiniteSplittings}
\begin{align}
\Z/n&=\prod_{p}\Z/p^{\ord_p n},
\label{eqn:ProfiniteSplittings-Finite}\\
\hat{\Z}&=\prod_{p}\Z_{p}.
\label{eqn:ProfiniteSplittings-Infinite}
\end{align}
\end{subequations}
Applying the Transitivity Sequence and using that $R_{(p)}$ is \'etale over 
$R$, we obtain that at each prime $p$
$$ \HG^*(R_{(p)}| \Z_{(p)}; \hat{\Z}) \cong \HG^*(R | \Z; \hat{\Z}).$$
Therefore by
Corollary~\ref{prop:ContGH-MilnorSeq} we have
\[
\HG^*(R | \Z;\hat{\Z})=\prod_{p}\HG^*(R|\Z;\Z_p)
\cong \prod_{p}\HG^*(R_{p}|\Z_{p};\Z_{p}).
\]
For the second isomorphism, using Corollary~\ref{prop:ContGH-MilnorSeq}
and the  linear compactness of $\Gamma$-cohomology provided by
Proposition~\ref{prop:Cohom-LinCompact}, we can express $\HG^*(R|\Z;\Z_p)$
as the inverse limit of the groups $\HG^*(R|\Z;\Z/p^n)$. Here the coefficients
$\Z/p^n$ eliminate the effect of all the $p$-divisible elements, therefore
$\HG^*(R|\Z;\Z/p^n)$ reduces to $\HG^*(R_p|\Z_p;\Z/p^n)$, where $R_p$
denotes the $p$-adic completion of $R$.

Now for each $k\geq1$, the assumption that $R/p^k$ is formally \'etale over
$\Z/p^k$ ensures that
\[
\HG^*(R/p^k|\Z/p^k;\Z/p^k)=0.
\]
Therefore we obtain
\[
\HG^*(R_p|\Z_{p};\Z_{p})=
\lim_{k}\HG^*(R/p^k|\Z/p^k;\Z/p^k)=0
={\lim_{k}}^1\;\HG^*(R/p^k|\Z/p^k;\Z/p^k).
\]
For each~$n$, Proposition~\ref{prop:LESk-k^} implies that
\begin{align*}
\HG^n(R_{(p)}|\Z_{(p)};\Z_{(p)})
&=\HG^{n-1}(R_{(p)}|\Z_{(p)};\Z_p/\Z_{(p)}), \\
\HG^n(R|\Z;\Z)&=\HG^{n-1}(R|\Z;\hat{\Z}/\Z).
\end{align*}
As $\Z_p/\Z_{(p)}$ and $\hat{\Z}/\Z$ are $\Q$-vector spaces, for all $n\neq0$
we obtain
\begin{align*}
\HG^n(R_{(p)}|\Z_{(p)};\Z_p/\Z_{(p)})&\iso\HG^n(R_\Q|\Q;\Z_p/\Z_{(p)}) \\
\intertext{and similarly}
\HG^n(R|\Z;\hat{\Z}/\Z)&\iso\HG^n(R_\Q|\Q;\hat{\Z}/\Z).
\qedhere
\end{align*}
\end{proof}
\begin{cor}\label{cor:HG-NumPolys}
We have
\[
\HG^n(\As|\Z;\Z)=\HG^n(\A|\Z;\Z)=
\begin{cases}
\hat{\Z}/\Z& \text{\rm if $n=1$}, \\
0& \text{\rm otherwise}.
\end{cases}
\]
For each prime $p$,
\[
\HG^n(\As_{(p)}|\Z_{(p)};\Z_{(p)})=\HG^n(\A_{(p)}|\Z_{(p)};\Z_{(p)})=
\begin{cases}
{\Z}_p/\Z_{(p)}& \text{\rm if $n=1$}, \\
0& \text{\rm otherwise}.
\end{cases}
\]
\end{cor}
\begin{proof}
Since $\A[w^{-1}]$ is \'etale over $\A$, both of the $\Gamma$-cohomology
groups $\HG^*(\A[w^{-1}]|\A;\Z)$ and $\HG^*(\A_{(p)}[w^{-1}]|\A_{(p)};\Z)$
vanish. The Transitivity Theorem~\cite[3.4]{RoWh} implies that there are
isomorphisms
\begin{align*}
\HG^*(\A|\Z;\Z)&\iso\HG^*(\A[w^{-1}]|\Z;\Z)=\HG^*(\As|\Z;\Z), \\
\HG^*(\A_{(p)}|\Z_{(p)};\Z_{(p)})
&\iso\HG^*(\A_{(p)}[w^{-1}]|\Z_{(p)};\Z_{(p)})
=\HG^*(\As_{(p)}|\Z_{(p)};\Z_{(p)}),
\end{align*}
hence it suffices to prove the result for $\A$ and $\A_{(p)}$.

Corollary~\ref{cor:NumPolys-Local} ensures that $\A/p^k\A$ and therefore
$\A_{(p)}/p^k\A_{(p)}$ as well is formally \'etale over $\Z/p^k$ for all
$k\geq1$. Now Corollaries~\ref{cor:NumPolys-Local}(d) and 
\ref{cor:Hom-LinCompact} together 
guarantee that the cochains for $\Gamma$-cohomology fulfil the linear
compactness requirements of Theorem~\ref{thm:Z-algebras}. Thus we can apply 
this theorem and obtain the vanishing
result for $\Gamma$-cohomology in dimensions different from one.

By~\cite[theorem 4.1]{BR&Ro} and the fact that $\Z_p/\Z_{(p)}$ and 
$\hat{\Z}/\Z$ are $\Q$-vector spaces, we have
\[
\HG^*(\A_{(p)}|\Z_{(p)};\Z_p/\Z_{(p)})
=\HG^*(\A\otimes\Q|\Q;\Z_p/\Z_{(p)})
=\HG^*(\Q[w]|\Q;\Z_p/\Z_{(p)})=\Z_p/\Z_{(p)}
\]
and
\[
\HG^*(\A|\Z;\hat{\Q}/\Q)=\HG^*(\A\otimes\Q|\Q;\hat{\Z}/\Z)
=\HG^*(\Q[w]|\Q;\hat{\Z}/\Z)=\hat{\Z}/\Z.
\]
Thus we obtain
\begin{align*}
\HG^n(\A_{(p)}|\Z_{(p)};\Z_{(p)})&=
\begin{cases}
\Z_p/\Z_{(p)}& \text{if $n=1$}, \\
0& \text{otherwise}, \\
\end{cases} \\
\intertext{and}
\HG^n(\A|\Z;\Z)&=\begin{cases}
\hat{\Z}/\Z& \text{if $n=1$}, \\
0& \text{otherwise},
\end{cases}
\end{align*}
as claimed.
\end{proof}

\begin{rem}\label{thm:AQ-A&As}
Notice, that for the calculations of $\Gamma$-cohomology above we
used the formal properties of $\Gamma$-cohomology. As Andr\'e-Quillen
cohomology satisfies analoguous properties, we can transfer the above
results to obtain the following:
\begin{align*}
\AQ_n(\As|\Z;\Z)=\AQ_n(\A|\Z;\Z)&=
\begin{cases}
\Q&\text{\rm if $n=0$},\\
0&\text{\rm if $n\neq0$},
\end{cases} \\
\AQ^n(\As|\Z;\Z)=\AQ^n(\A|\Z;\Z)&=
\begin{cases}
\hat{\Q}/\Q&\text{\rm if $n=1$},\\
0&\text{\rm if $n\neq1$}.
\end{cases}
\end{align*}
\end{rem}

The results from Section~\ref{sec:Z/(p-1)invariants} allow us to
calculate the $\Gamma$-cohomology of ${}^\zeta\As_{(p)}$ over $\Z_{(p)}$
directly as was done above for $\As$. Alternatively, we may use the
fact that the extension $\As_{(p)}/{}^\zeta\As_{(p)}$ is \'etale since
it has the form $B/A$, where $B=A[t]/(t^{p-1}-v)$ for a unit $v\in A$,
where $A$ is a $\Z_{(p)}$-algebra. We can now determine the
$\Gamma$-cohomology of ${}^\zeta\As_{(p)}$ since the Transitivity Theorem
of~\cite[3.4]{RoWh} gives
\begin{prop}\label{prop:HG-omegaAs(p)}
For an odd prime $p$,
\[
\HG^*({}^\zeta\As_{(p)}|\Z_{(p)};\Z_{(p)})=
\HG^*(\As_{(p)}|\Z_{(p)};\Z_{(p)}).
\]
\end{prop}

\section{Applications to $E_\infty$ structures on $K$-theory}
\label{sec:EinftyStructures}

Robinson~\cite{Ro:Einfty} has developed an obstruction theory for
$E_\infty$ structures on a homotopy commutative ring spectrum~$E$.
Provided $E$ satisfies the following form of the  K\"unneth and
universal coefficient
theorems for $E_*E$ 
$$ E^*(E^{\wedge n}) \cong \Hom_{E_*}(E_*E^{\otimes n}, E_*), $$
then the obstructions lie in groups
\[
\HG^{n,2-n}(E_*E| E_*;E_*),
\]
while the extensions are determined by classes in
\[
\HG^{n,1-n}(E_*E| E_*;E_*).
\]
Here the bigrading $(s,t)$ involves cohomological degree~$s$ and
internal degree~$t$. Moreover, the relevant values of $n$ are for
$n\geq 3$.

We want to apply this to the cases of complex $KU$-theory and the
Adams summand $E(1)$ of $KU_{(p)}$ at a prime~$p$. Recall that
\[
KU_*=\Z[t,t^{-1}],
\quad
{KU_{(p)}}_*=\Z_{(p)}[t,t^{-1}],
\quad
E(1)_*=\Z_{(p)}[u,u^{-1}],
\]
where $t\in KU_2$ and $u\in E(1)_{2(p-1)}$. Our next result implies
that the relevant conditions mentioned above are both satisfied for
$KU$ and $E(1)$.
\begin{prop}\label{prop:KU*KU}
There are isomorphisms of rings {\rm(}in fact, of Hopf algebras{\rm)}
\[
KU_0KU\iso\As,
\quad
{KU_{(p)}}_0KU_{(p)}\iso\As_{(p)},
\quad
E(1)_0E(1)\iso{}^\zeta\As_{(p)}.
\]
Hence,
\[
KU_*KU\iso KU_*\otimes\As,\quad {KU_{(p)}}_*KU_{(p)}\iso {KU_{(p)}}_*
\otimes \As_{(p)},\quad  E(1)_*E(1)\iso E(1)_*\otimes{}^\zeta\As_{(p)}.
\]
\end{prop}
\begin{proof}
The isomorphisms for $KU$ and $KU_{(p)}$ can be found in~\cite[p.~392]{BCRS}.

Consider $E(1)_*E(1)$, the algebra of cooperations for $E(1)$. Since
$E(1)$ is Landweber exact, we have
\[
E(1)_*E(1)\cong E(1)_*\otimes_{BP_*}BP_*BP\otimes_{BP_*}E(1)_*
\cong E(1)_*[t_1,t_2,\ldots,V_1,V_1^{-1}]/(\text{relations}),
\]
where $V_1$ denotes the right unit $\eta_r$ applied to $v_1$ and the
variables $t_i$ stem from $BP_*BP$. We also write $\bar{w}=v_1^{-1}V_1$.
The relation $v_1+pt_1-V_1=0$ in $E(1)_*E(1)$ gives rise to
\[
1-v_1^{-1}V_1=-pv_1^{-1}t_1\in E(1)_0E(1),
\]
hence on setting $\bar{e}_1=v_1^{-1}t_1$ we have
\[
\bar{w}-1=p\bar{e}_1.
\]
Now we may inductively define
\[
\bar{e}_m=v_1^{-p^{m-1}-\cdots-p-1}t_m .
\]
The higher relations
\[
v_1t_k^{p}-v_1^{p^k}t_k+pt_{k+1}=0\quad(k\geq1)
\]
can be used to prove the desired relations for the $\bar{e}_m$.
Taking the $p$-th power we have
\[
\bar{e}_m^p=v_1^{-p^m-\cdots-p}t_m^p.
\]
Multiplying the relation
\[
v_1t_m^p-v_1^{p^m}t_m=pt_{m+1}
\]
by $v_1^{-p^m-\cdots-p-1}$, we obtain
\[
\bar{e}_m^p-v_1^{p^m-p^m-p^{m-1}-\cdots-p-1}t_m
=pv_1^{-p^m-\cdots-p-1}t_{m+1}
\]
which is precisely
\[
\bar{e}_m^p-\bar{e}_m=p\bar{e}_{m+1}.
\qedhere
\]
\end{proof}\noindent
Flat base-change leads to isomorphisms
\begin{align*}
\HG^{*,*}(KU_*KU| KU_*;KU_*)&\cong\HG^*(\As|\Z;\Z)\otimes KU_*, \\
\HG^{*,*}({KU_{(p)}}_*{KU_{(p)}}|{KU_{(p)}}_*;{KU_{(p)}}_*)
     &\cong\HG^*(\As_{(p)}|\Z_{(p)};\Z_{(p)})\otimes{KU_{(p)}}_*, \\
\HG^{*,*}(E(1)_*E(1)| E(1)_*;E(1)_*)
         &\cong\HG^*({}^\zeta\As_{(p)}|\Z_{(p)};\Z_{(p)})\otimes E(1)_*.
\end{align*}
With the help of Corollary~\ref{cor:HG-NumPolys} we can therefore deduce
the following.
\begin{thm}\label{thm:K-Obstructions}
For a prime $p$ and $n\geq 2$,
\begin{align*}
\HG^{n,2-n}(KU_*KU| KU_*;KU_*)
&=0=
\HG^{n,1-n}(KU_*KU| KU_*;KU_*), \\
\HG^{n,2-n}({KU_{(p)}}_*KU_{(p)}|{KU_{(p)}}_*;{KU_{(p)}}_*)
&=0=
\HG^{n,1-n}({KU_{(p)}}_*KU_{(p)}|{KU_{(p)}}_*;{KU_{(p)}}_*), \\
\HG^{n,2-n}(E(1)_*E(1)| E(1)_*;E(1)_*)
&=0=
\HG^{n,1-n}(E(1)_*E(1)| E(1)_*;E(1)_*).
\end{align*}
Hence $KU$, $KU_{(p)}$, and $E(1)$ each have a unique $E_\infty$
structure.
\end{thm}

It is a rather old question whether the connective Adams summand,
often denoted by~$\ell$, is an $E_\infty$ spectrum. The $E_\infty$
ring spectrum machinery developed in~\cite{May:EinftyRings} yields
the following general result.
\begin{thm}\label{thm:ConnCov}
For any $E_\infty$ ring spectrum $E$, the connective cover $e\lra E$
possesses a model as an $E_\infty$ ring spectrum.
\end{thm}
\begin{proof}
Proceeding as in~\cite[prop.~VII.3.2]{May:EinftyRings}, we first take
the underlying zeroth space $E_0$ of the $E_\infty$ ring spectrum~$E$,
then build a prespectrum $T(E_0)$ out of it using a bar construction
which consists of suspensions and the monad for the little convex body
(partial) operad. Finally we apply the spectrification functor (there
called $\Omega^{\infty}$) to $T(E_0)$. By~\cite[prop.~VII.3.2]{May:EinftyRings},
this has the correct homotopy groups and is an $E_\infty$ ring spectrum.
\end{proof}

Applying this result, we obtain a canonical $E_\infty$ model for the
connective cover $\ell\lra E(1)$.
\begin{prop}\label{prop:AdamsSummand}
There is at least one $E_\infty$ structure on the connective Adams
summand~$\ell$.
\end{prop}
\begin{rem}\label{rem:AdamsSummand-padic}
After $p$-completion, we obtain an $E_\infty$ structure on the
$p$-completed connective Adams summand $\ell\sphat_p$\,. In subsequent
work we have shown  that this $E_\infty$ structure coincides with the
one  constructed by McClure and Staffeldt in~\cite{McC&S} using
algebraic $K$-theory.
\end{rem}

\section{$E_\infty$ structures on $KO$}\label{sec:KO}

The case of $KO$ can be treated by similar methods but involves
somewhat more delicate considerations because of the presence of
$2$-torsion in $KO_*$. Recall that
\begin{equation}\label{eqn:KO-CoeffRing}
KO_*=\Z[h,y,w,w^{-1}]/(2h,h^3,hy,y^2-4w),
\end{equation}
where $h\in KO_1$, $y\in KO_4$ and $w\in KO_8$. We will also require
the graded $\Q$-vector space $V_*=KO_*\otimes\hat{\Z}/\Z$.

We will prove the following algebraic result.
\begin{thm}\label{thm:HGamma-KO}{\ }
\begin{enumerate}
\item[(a)]
For any prime $p$ and $k\geq1$, we have
\[
\HG^*(KO_0KO/p^k|\Z/p^k;\Z/p^k)=0,
\]
and
\[
\cHG^*(KO_0KO\sphat_p|\Z_p;\Z_p)=\HG^*(KO_0KO\sphat_p|\Z_p;\Z_p)=0.
\]
\item[(b)]
We have
\[
\HG^{n,*}(KO_*KO| KO_*;KO_*)=
\begin{cases}
V_*& \text{\rm if $n=1$}, \\
0& \text{\rm otherwise}.
\end{cases}
\]
\end{enumerate}
\end{thm}

Using this, the obstruction theory $E_\infty$ structures and localization
yield our result on $E_\infty$ structures for $KO$.
\begin{thm}\label{thm:KO-Einfty}
$KO$\/ and for each prime $p$, $KO_{(p)}$\/ and $KO\sphat_p$\/ all have
unique $E_\infty$ structures.
\end{thm}

To prove Theorem~\ref{thm:HGamma-KO}, we begin with a composite result
distilled from~\cite{AHS} and~\cite[p.~162]{JFA:InfLoopSpaces}.
\begin{thm}\label{thm:KO*KO}{\ }
\begin{enumerate}
\item[(a)]
$KO_*KO$ is a free $KO_*$-module on countably many generators lying
in $KO_0KO$.
\item[(b)]
The natural homomorphism
\[
KO_0KO\lra KU_0KU\xrightarrow{\iso}\As
\]
is a split monomorphism whose image is
\[
\{f(w)\in\As:f(-w)=f(w)\}\subset\As.
\]
\item[(c)]
For each prime $p$ and $k\geq1$,
\[
KO_0KO/p^k\lra KU_0KU/p^k\iso\As/p^k
\]
is a split monomorphism of $\Z/p^k$-modules.
\end{enumerate}
\end{thm}

Consider the short exact sequence
\[
0\ra KO_*\lra KO_*\otimes\hat{\Z}\lra KO_*\otimes\hat{\Z}/\Z\ra0.
\]
Since $\hat{\Z}/\Z$ is a $\Q$-vector space and we have the splitting
of~\eqref{eqn:ProfiniteSplittings-Infinite}, we can reformulate the
above exact sequence to obtain that the sequence
\begin{equation}\label{eqn:KO-ExactSeq}
0\ra KO_*\lra\prod_p KO_*\otimes\Z_p\lra V_*\ra0,
\end{equation}
is exact. Here $V_*$ is defined above and we have used the fact that
each group $KO_n$ is finitely-generated. The application of
$\Gamma$-cohomology of $KO_*KO$ to this sequence yields a long exact
sequence which relates $\HG^{*,*}(KO_*KO| KO_*;KO_*)$ to
$\HG^{*,*}(KO_*KO| KO_*;\prod_p KO_*\otimes\Z_p)$ and
$\HG^{*,*}(KO_*KO| KO_*;V_*)$.

Now, for the part with coefficients in $V_*$ we have
\[
\HG^{*,*}(KO_*KO| KO_*;V_*)\iso
\HG^{*,*}(KO_*KO\otimes\Q| KO_*\otimes\Q;V_*),
\]
and
\[
KO_*\otimes\Q=\Q[y,y^{-1}],
\quad
KO_*KO\otimes\Q=\Q[y,y^{-1},z,z^{-1}].
\]
By~\cite{BR&Ro}, as in Corollary~\ref{cor:HG-NumPolys} we find that
\begin{equation}\label{eqn:HG-V*}
\HG^{n,*}(KO_*KO| KO_*;V_*)=
\begin{cases}
V_*& \text{if $n=1$}, \\
0& \text{otherwise}.
\end{cases}
\end{equation}
For a fixed algebra, the cochain complex for $\Gamma$-cohomology
commutes with limits taken over the coefficient module, therefore
$\Gamma$-cohomology commutes with products of coefficient modules
and the splitting $\hat{\Z}=\prod_p\Z_p$ leads to
\[
\HG^{*,*}(KO_*KO| KO_*;KO_*\otimes\hat{\Z})
\cong \prod_p\HG^{*,*}(KO_*KO| KO_*;KO_*\otimes\Z_p).
\]
For each prime $p$ we obtain a short exact sequence,
\begin{multline*}
0\ra{\lim_k}^1\; \HG^{*-1,*}(KO_*KO| KO_*;KO_*/p^k)
\lra\HG^{*,*}(KO_*KO| KO_*;KO_*\otimes\Z_p) \\
\lra\lim_k\HG^{*,*}(KO_*KO| KO_*;KO_*/p^k)\ra0.
\end{multline*}
When $p>2$, we are reduced to considering
\[
\HG^{*,*}(KO_*KO| KO_*;KO_*/p^k)
=\HG^{*,*}(KO_*KO/p^k| KO_*/p^k;KO_*/p^k),
\]
which can be determined by the methods of Section~\ref{sec:Z/(p-1)invariants}
using the subgroup $\{\pm1\}\leq\Z_p^\times$ in place of the group of all
$(p-1)$-st roots of unity. The result is that
\[
\HG^{*,*}(KO_*KO| KO_*;KO_*/p^k)=0,
\]
whence
\[
\HG^{*,*}(KO_*KO| KO_*;KO_*\otimes\Z_p)=0.
\]

The case $p=2$ requires a more intricate analysis. First we identify 
$KO_0KO_{(2)}$ and the quotients $KO_0KO/2^k$ as rings of functions.
\begin{thm}\label{thm:KO0KO-Fns}{\ }
\begin{enumerate}
\item[(a)]
There is an isomorphism of rings
\[
KO_0KO_{(2)}\iso
\{f(w)\in\Q[w,w^{-1}]:f\Z_{(2)}^\times\subset\Z_{(2)},\;f(-w)=f(w)\}
\subset\As_{(2)}.
\]
\item[(b)]
For each $k\geq1$, there is an isomorphism of rings
\[
KO_0KO/2^k\iso\Cont(1+8\Z_2,\Z/2^k),
\]
where $\Cont(1+8\Z_2,\Z/2^k)$ denotes the space of continuous maps
from $1+8\Z_2\subset\Z_2^\times\subset\Z_2$ with its $2$-adic topology
to $\Z/2^k$ with the discrete topology.
\item[(c)]
There is an isomorphism of rings
\[
KO_0KO\sphat_2\iso\Cont(1+8\Z_2,\Z_2),
\]
the space of continuous maps from $1+8\Z_2$ to $\Z_2$.
\item[(d)]
The algebras $(KO_0KO/2^k,\Z/2^k)$ and $(KO_0KO\sphat_2,\Z_2)$ are
formally \'etale.
\end{enumerate}
\end{thm}
\begin{proof}
The methods of~\cite{AB:p-adic} apply here, and we leave verification
of the details to the reader.

The squaring map $\Z_2^\times\lra\Z_2^\times$ has image $1+8\Z_2$, hence
a polynomial $f(w)\in\As_{(2)}$ satisfying $f(-w)=f(w)$ corresponds to
a continuous function $1+8\Z_2\lra\Z_2$. By compactness of the domain,
$\Cont(1+8\Z_2,\Z/2^k)$ consists of locally constant functions. If we
express $x\in\Z_2$ in the form
\[
x=x_0+x_12+x_22^2+\cdots+x_n2^n+\cdots,
\]
where $x_i=0,1$, then the functions
\[
\xi_i\:1+8\Z_2\lra\Z_2;\quad\xi_i(x)=x_i
\]
are locally constant and give rise to $\Z/2^k$-algebra generators of
$\Cont(1+8\Z_2,\Z/2^k)$. They also satisfy the relations
\[
\xi_i^2=\xi_i,
\]
and the distinct monomials
\[
\xi_0^{r_0}\xi_1^{r_1}\dotsm\xi_d^{r_d}
\quad(r_i=0,1)
\]
form a $\Z/2^k$-basis. This implies that the $\Z/2^k$-algebra
$\Cont(1+8\Z_2,\Z/2^k)$ is formally \'etale. Similar considerations
apply to the topological algebra $KO_0KO\sphat_2$.
\end{proof}

Collecting together the results of the above discussion (in particular
Theorem~\ref{thm:KO0KO-Fns}(d)) we obtain the case $p=2$ of
Theorem~\ref{thm:HGamma-KO}(a). The proof of Theorem~\ref{thm:HGamma-KO}(b)
makes use of the long exact sequence resulting from~\eqref{eqn:KO-ExactSeq}
and~\eqref{eqn:HG-V*}.

We remark that rather than working modulo powers of~$2$, it is
also possible to consider powers of the maximal ideal $(2,h,y)\ideal KO_*$
and then we obtain
\begin{prop}\label{prop:KO*KO-HG-(2,y)}
For $k\geq1$,
\[
\HG^*(KO_*KO/(2,h,y)^k| KO_*/(2,h,y)^k;KO_*/(2,h,y)^k)=0.
\]
and
\begin{multline*}
\cHG^{*,*}(KO_*KO\sphat_{(2,h,y)}|(KO_*)\sphat_{(2,h,y)};
(KO_*)\sphat_{(2,h,y)}) \\
=\HG^{*,*}(KO_*KO\sphat_{(2,h,y)}|(KO_*)\sphat_{(2,h,y)};
(KO_*)\sphat_{(2,h,y)})
=0.
\end{multline*}
\end{prop}

\section{$E_\infty$ structures on the $I_n$-adic completion of $E(n)$}
\label{sec:E(n)-InadicComp}

In this section we describe what we can prove about $E_\infty$ structures
on the $I_n$-adic completion of Johnson-Wilson spectrum $E(n)$ for a prime~$p$
and $n\geq1$.

The coefficient ring
\[
E(n)_*=\Z_{(p)}[v_1,\ldots,v_{n-1},v_n,v_n^{-1}]
\]
is Noetherian and contains the maximal ideal
\[
I_n=(p,v_1,\ldots,v_{n-1})\ideal E(n)_*.
\]
Here the $v_i$ denote the images of the Araki generators of $BP_*$ and we
sometimes write $v_0=p$. There is a commutative ring spectrum $\hat{E(n)}$
for which the coefficient ring $\hat{E(n)}_*$ is the $I_n$-adic completion
of $E(n)_*$, \ie, its completion at $I_n$. It is known from~\cite{AB&UW,HS}
that $\hat{E(n)}$ is the $K(n)$-localization of $E(n)$. We also know
from~\cite{AB:Ainfty} that for each prime~$p$, $\hat{E(n)}$ possesses a
unique $A_\infty$ structure and the canonical map
$\hat{E(n)}\lra\hat{E(n)}/I_n\homeq K(n)$ to the $n$-th Morava $K$-theory
is a map of $A_\infty$ ring spectra for any of the $A_\infty$ structures
on $K(n)$ shown to exist in~\cite{Ro:Ainfty}. Actually these results were
only claimed for odd primes but the arguments also work for the prime~$2$.

\begin{prop}
Possible obstructions for an $E_\infty$ structure on the completed Johnson-Wilson
spectra $\hat{E(n)}$ live in the continuous $\Gamma$-cohomology groups
$$ \cHG^{*,*}(\hat{E(n)}_*\hat{E(n)}|\hat{E(n)}_*;\hat{E(n)}_*).$$
\end{prop}
\begin{proof}
For $\hat{E(n)}$ we have a
\emph{continuous} universal coefficient theorem, \ie, possible obstructions 
live in the \emph{continuous} $\hat{E(n)}$-cohomology of
$(X_m)_+\ltimes_{\Sigma_m}E^{\wedge m}$, where $X_m$ is a filtration
quotient of an $E_\infty$ operad as described in~\cite[section~5.1]{Ro:Einfty}.
These cohomology groups can be identified with the continuous
$\hat{E(n)}_*$-homomorphisms from the corresponding $\hat{E(n)}$-homology
groups (compare~\cite[proposition 5.4]{Ro:Einfty}
and~\cite[\S1]{AB:Ainfty}).  This proves
the claim. 
\end{proof}
For each $\ell\geq0$, Proposition~\ref{prop:ContGH-fgmodule} yields a short 
exact sequence
\begin{multline} \label{eqn:milnorhaten}
0\ra{\lim_k}^1\;\HG^{\ell-1,*}(E(n)_*E(n)/I_n^k| E(n)_*/I_n^k;E(n)_*/I_n^k) \\
\lra\cHG^{\ell,*}(\hat{E(n)}_*\hat{E(n)}|\hat{E(n)}_*;\hat{E(n)}_*) \\
\lra{\lim_k}\;\HG^{\ell,*}(E(n)_*E(n)/I_n^k| E(n)_*/I_n^k;E(n)_*/I_n^k)
\ra0.
\end{multline}
\begin{thm}\label{thm:E(n)Ink-Gammacohomology=0}
The $E(n)_*/I_n^k$-algebra $E(n)_*E(n)/I_n^k$ is formally \'etale. Hence
the $\Gamma$-cohomology of $E(n)_*E(n)/I_n^k$ over $E(n)_*/I_n^k$ is trivial,
\[
\HG^{*,*}(E(n)_*E(n)/I_n^k| E(n)_*/I_n^k;\hat{E(n)}_*/I_n^k)=0.
\]
\end{thm}
\begin{proof}
First we show that the algebra $E(n)_*E(n)/I_n^k$ is formally \'etale. In the
following we use the notation chosen in \cite{AB:Ainfty}.
As in the proof of ~\cite[lemma 3.4]{AB:Ainfty}, we can apply the 
infinite-dimensional Hensel lemma
(see the proof of our Corollary~\ref{cor:NumPolys-Local}) to split
$E(n)_*E(n)/I_n^k$ into an infinite tensor product of $E(n)_*/I_n^k$-algebras,
\[
E(n)_*E(n)/I_n^k=
\bigotimes_{j\geq1}
E(n)_*/I_n^k[S_j]/(v_nS_j^{p^n} - v_n^{p^j}S_j).
\]
We can write $E(n)_*E(n)/I_n^k$ as a colimit of finite tensor products,
\[
E(n)_*E(n)/I_n^k=\colim_m\bigotimes_{j=1}^m
E(n)_*/I_n^k[S_j]/(v_nS_j^{p^n}-v_n^{p^j}S_j).
\]
We claim that each algebra $E(n)_*/I_n^k[S_j]/(v_nS_j^{p^n}-v_n^{p^j}S_j)$
is \'etale over $E(n)_*/I_n^k$. Notice that it is flat over $E(n)_*/I_n^k$
and is finitely-generated by $S_j$. As the ground ring $E(n)_*/I_n^k$ is
Noetherian, the only thing that remains to be shown is that the module
of K\"ahler differentials is trivial.

The K\"ahler differentials are generated by the symbol $dS_j$, but in
$E(n)_*E(n)/I_n^k$ we have the relation $v_nS_j^{p^n}=v_n^{p^j}S_j$. The
residue class of the element $v_n\in E(n)_*$ is a unit in the ring
$E(n)_*/I_n$ and thus we can deduce
\[
dS_j=v_n^{1-p^j}d(S_j^{p^n})=p^nv_n^{1-p^j}S_j^{p^n-1}dS_j.
\]
Iteration of this relation $t$ times, where $t$ is an integer such
that $tn \geq k$, implies that $dS_j$ is zero, since in the quotient
$E(n)_*/I_n^k$, $p^k$ is zero.

Now by Lemma~\ref{lem:HG-Formallyetale} $\Gamma$-homology commutes with
colimits, therefore
\[
\HG_{*,*}(E(n)_*E(n)/I_n^k| E(n)_*/I_n^k;\hat{E(n)}_*/I_n^k)=0
=\HG^{*,*}(E(n)_*E(n)/I_n^k| E(n)_*/I_n^k;\hat{E(n)}_*/I_n^k).
\]
This completes the proof of Theorem~\ref{thm:E(n)Ink-Gammacohomology=0}.
\end{proof}

Using (\ref{eqn:milnorhaten}) and the fact that the completion of 
$\hat{E(n)}_*\hat{E(n)}$ is free on a countable basis 
\cite[theorem 1.1]{AB:left}, \cite{AB:p-adic}, we obtain
\begin{thm}\label{thm:hatE(n)-Einftystructure}
For $p$ a prime and $n\geq1$, the spectrum $\hat{E(n)}$ possesses
a unique $E_\infty$ structure.
\end{thm}

Using the ideas of Section~\ref{sec:GammCoh}, we can also deduce
\begin{thm}\label{thm:AQ-E(n)E(n)^}
For $n\geq1$ and $k\geq1$, we have
\begin{align*}
\AQ_*(E(n)_*E(n)/I_n^k| E(n)_*/I_n^k;E(n)_*/I_n^k)
&=0=\AQ^*(E(n)_*E(n)/I_n^k| E(n)_*/I_n^k;E(n)_*/I_n^k), \\
\AQ_*(E(n)_*E(n)\sphat_{I_n}|\hat{E(n)}_*;\hat{E(n)}_*)
&=0=\AQ^*(E(n)_*E(n)\sphat_{I_n}|\hat{E(n)}_*;\hat{E(n)}_*).
\end{align*}
\end{thm}


\begin{rem}\label{rem:E(n)-Einftystructure}
Extending Theorem~\ref{thm:hatE(n)-Einftystructure} to cover $E(n)$
for $n>1$ does not appear to be straightforward. The following two 
problems arise.
\begin{itemize}
\item
There is the question of whether $E(n)_*E(n)$ is a free $E(n)_*$-module
when $n>1$. If $E(n)$ does not have a universal coefficient theorem,
then the obstructions to building an $E_\infty$ structure on $E(n)$
would live in $E(n)$-cohomology which might not then be identifiable
with $\Gamma$-cohomology. In~\cite{AB:EnEn}, the first author showed
that the cooperation algebra of the $I_n$-localization of $E(n)$,
$E(n)_{I_n}$, is free over ${E(n)_*}_{I_n}$, so it does have a universal
coefficient theorem and the above problem is overcome.
\item
$\Gamma$-cohomology of $E(n)_*E(n)$ is non-trivial in positive degrees.
Even for $n=2$ there are polynomial generators in $E(2)_*E(2)$ which
lead to non-trivial classes in $\Gamma$-cohomology.
\end{itemize}
We aim to return to the existence of $E_\infty$ structures on $E(n)$
and $E(n)_{I_n}$ in future work.
\end{rem}

We end this section with some remarks on suitably completed versions
of elliptic cohomology. Here $\Ell$ denotes the level~$1$ version of
elliptic cohomology of Landweber, Ravenel and Stong~\cite{LRS} and we
focus on the supersingular completions of~\cite{AB:homell}. Our above
techniques together with results from~\cite{AB:EllEll} yield the
following.
\begin{thm}\label{thm:EllCoh}
For each prime $p>3$, the supersingular completions $\Ell\sphat_{(p,E_{p-1})}$
and $\Ell\sphat_{\mathcal{P}}$ for each maximal ideal $\mathcal{P}\ideal(\Ell_*)_{(p)}$
containing $(p,E_{p-1})$, have unique $E_\infty$ structures.
\end{thm}

An analogous result applies to the $K(1)$-localization of $\Ell$ studied
in~\cite{AB:p-hecke} and more recently by M.~Hopkins.

\section{An obstruction theory for the coherence of maps}\label{sec:CohMaps}

For the following, we  need to work in a good category of spectra with
a symmetric monoidal smash product, for example that of~\cite{EKMM}.
Where necessary, all ring spectra will be assumed to be fibrant.

Let $E$ and $F$ be two $E_\infty$ ring spectra over the $E_\infty$
operad~$\T$ from~\cite[section~5.1]{Ro:Einfty} and let $f\:E\lra F$
be a map of commutative ring spectra, \ie, the map $f$ commutes with
the multiplication maps $\mu_E$ and $\mu_F$ up to homotopy,
\[
\mu_F\circ f\wedge f\homeq f\circ\mu_E,
\]
and similar coherence properties exist with respect to the homotopies
for associativity and commutativity on $E$ and $F$. The aim of the
following discussion is to give criteria, when the map $f$ can be made
into a map which is compatible with the $\T$-algebra structures on $E$
and $F$ up to homotopy. For $A_\infty$ structures the analogous question
was addressed in~\cite{Ro:Ainfty}.

From now on we will use the notation of~\cite{Ro:Einfty}. The topological
operad $\T$ is filtered by subspaces $\nabla^m\T(n)\subset\T(n)$. Let
$\theta_E$ and $\theta_F$ be the action maps of the operad $\T$ on $E$
and $F$.

Consider the sequence of topological spaces
\[
\nabla^m\T (n)\hookrightarrow
\nabla^{m+1}\T(n)\lra
\nabla^{m+1}\T(n)/\nabla^m\T(n)\cup\partial\nabla^{m+1}\T(n),
\]
where $\partial\nabla^{m+1}\T(n)$ is the part of $\T(n)$ which
is determined by compositions in the operad of elements coming
from lower filtration degrees.
\begin{thm}\label{thm:CohExtn}{\ }
\begin{enumerate}
\item[(a)]
If $\HG^{n,2-n}(F_*E| F_*;F_*)=0$ for all $n\geq3$, then $f$
can be turned into a map that satisfies
\[
f\circ\theta_E\homeq\theta_F\circ\underbrace{f\wedge\ldots\wedge f}_m\:
\T(m)\ltimes_{\Sigma_m}E^{\wedge m}\lra F
\]
for all $m$.
\item[(b)]
If in addition $\Hom_{F_*}(F_*E,F_*)\cong\Hom_{E_*}(E_*E,F_*)$,
then it suffices to prove that
\[
\HG^{n,2-n}(E_*E| E_*;F_*)=0
\]
for all $n\geq3$.
\end{enumerate}
\end{thm}

The second condition is satisfied for instance if $F$ is projective
over~$E$, then $F_*E\cong F_*\otimes_{E_*}E_*E$ can be used to reduce
the module of $F_*$-linear morphisms to the module of $E_*$-linear
morphisms.
\begin{proof}
Assume $f$ satisfies the conditions up to filtration degree~$m$. In
order to extend $f$ coherently over the $(m+1)$-st filtration step,
we have to show that the condition of the theorem suffices to force~$f$
to fulfil
\[
f\circ\theta_E|\nabla^{m+1}\homeq
\theta_F|\nabla^{m+1}\circ\underbrace{f\wedge\ldots\wedge f}_n
\:\nabla^{m+1}\T(n)\ltimes_{\Sigma_n}E^{\wedge n}\lra F.
\]

The map $f\circ\theta_E|\nabla^{m+1}$ corresponds to an element
in $F^0(\nabla^{m+1}\T(n)\ltimes_{\Sigma_n}E^{\wedge n})$. Using
the long exact cohomology sequence corresponding to the sequence
of spaces
\[
\nabla^m\T(n)\cup
\partial\nabla^{m+1}\T (n)
\xrightarrow{i}
\nabla^{m+1}\T(n)
\xrightarrow{j}
\nabla^{m+1}\T(n)/\nabla^m\T(n)\cup\partial\nabla^{m+1}\T(n),
\]
we find that the difference element
\[
f\circ\theta_E|\nabla^{m+1}-\theta_F\circ f^{\wedge n}|\nabla^{m+1}
\]
maps to zero under $i^*$, thus it has to be in the image of $j^*$.
Consequently, if $j^*$ has trivial codomain, then this difference has
to be trivial as an element in $F$-cohomology. An argument showing that
the corresponding class in
\[
F^0(\nabla^{m+1}\T(n)/\nabla^m\T(n)\cup\partial\nabla^{m+1}\T(n)
\ltimes_{\Sigma_n}E^{\wedge n})
\]
has to be a cocycle in the complex for $\HG^*$ can be found
in~\cite{Ro:Einfty}. Thus if $\HG^{m,2-m}(F_*E| F_*;F_*)$ vanishes
in all degrees $m\geq3$, the potentially obstructing difference maps
$f\circ\theta_E-\theta_F\circ f^{\wedge n}$ have to be nullhomotopic.
\end{proof}

{}From the triviality of $\HG^n$ when $n > 1$ for complex $K$-theory and
its localization at a prime~$p$, we can deduce the following result.
\begin{thm}\label{thm:Kthy-ADamsOps}
For each $k$ integer prime to $p$, the $k$-th Adams operation
$\psi^k\:KU_{(p)}\lra KU_{(p)}$ can be refined to a coherent map with
respect to the $E_\infty$ structure given by the operad action of $\T$
on $KU_{(p)}$.
\end{thm}
\begin{proof}
The action of such an Adams operation $\psi^k$ on
${KU_{(p)}}_{2n}\iso\Z_{(p)}$ is given by multiplication by $k^n$,
thus it induces a different ${KU_{(p)}}_*KU_{(p)}$-module structure
on ${KU_{(p)}}_*$. This corresponds to taking $E=KU_{(p)}=F$ and the
map $\psi^k\:E\lra F$, then applying Theorem~\ref{thm:CohExtn}(b)
and using the fact that the relevant $\Gamma$-cohomology groups vanish,
this being a generalization of Theorem~\ref{thm:K-Obstructions} which
is proved in a similar way (this result depends crucially on the
vanishing of $\Gamma$-cohomology for formally \'etale extensions).
\end{proof}

Finally we have a result on the inclusion $j\:E(1)\lra KU_{(p)}$ of
the Adams summand into $p$-local $K$-theory which is a map of ring
spectra.
\begin{prop}\label{prop:E(1)->KUp}
$j$ gives rise to a coherent map of $E_\infty$ spectra.
\end{prop}
\begin{proof}
Using the Conner-Floyd isomorphism and the Landweber exactness of $E(1)$, 
the above argument can be adapted to prove that the relevant part of 
$\HG^{*,*}({KU_{(p)}}_*(E(1))|{KU_{(p)}}_*;{KU_{(p)}}_*)$ vanishes.
\end{proof}
\begin{rem}
With the aid of more machinery one can actually take the above arguments
to obtain the existence of \emph{strict} maps of $E_\infty$ ring spectra.
Using a comparison result of Basterra and the second
author~\cite[theorem 2.6]{BR},
we can identify $\Gamma$-cohomology groups with the obstruction groups
arising in the work of Goerss and Hopkins~\cite{GH}. Now the Goerss-Hopkins
obstruction theory~\cite[\S4]{GH} tells us that the vanishing of the
$\Gamma$-cohomology groups
$\HG^{*,*}({KU_{(p)}}_*(E(1))|{KU_{(p)}}_*;{KU_{(p)}}_*)$ and
$\HG^{*,*}({KU_{(p)}}_*(KU_{(p)})|{KU_{(p)}}_*;{KU_{(p)}}_*)$ implies that 
the Adams operations and the map $j$ give rise to maps of $E_\infty$ ring 
spectra.
\end{rem}


\begin{thebibliography}{12}\frenchspacing
\bibitem{JFA:InfLoopSpaces}
J.~F.~Adams,
Infinite Loop Spaces,
Annals of Mathematics Studies {\bf90},
Princeton University Press (1978).
\bibitem{Ad-Cl}
J.~F.~Adams \&  F.~W.~Clarke,
Stable operations on complex $K$-theory,
Ill. J. Math. {\bf21} (1977), 826--829.
\bibitem{AHS}
J.~F.~Adams, A.~S.~Harris \& R.~M.~Switzer,
Hopf algebras of co-operations for real and complex $K$-theory,
Proc. Lond. Math. Soc. {\bf23} (1971), 385--408.
\bibitem{AB:p-adic}
A.~Baker,
$p$-adic continuous functions on rings of integers,
J. Lond. Math. Soc. {\bf33} (1986), 414--20.
\bibitem{AB:p-hecke}
\bysame,
Elliptic cohomology, $p$-adic modular forms and Atkin's
operator $\mathrm{U}_p$, Algebraic topology, Proc. Int. Conf.,
Evanston/IL 1988, M.~Mahowald, S.~Priddy (eds), 
Contemp. Math. {\bf96} (1989), 33--38.
\bibitem{AB:homell}
\bysame,
On the homotopy type of the spectrum representing elliptic
cohomology,
Proc. Amer. Math. Soc. {\bf107} (1989), 537--48.
\bibitem{AB:Ainfty}
\bysame,
$A_\infty$ structures on some spectra related to Morava $K$-theory,
Quart. J. Math. Oxf. (2), {\bf42} (1991), 403--419.
\bibitem{AB:EllEll}
\bysame,
Operations and cooperations in elliptic cohomology,
Part I: Generalized modular forms and the cooperation algebra,
New York J. Math. {\bf1} (1995), 39--74.
\bibitem{AB:left}
\bysame,
A version of Landweber's filtration theorem for $v_n$-periodic Hopf
algebroids,  
Osaka J. Math. {\bf 32}  (1995),  no. 3, 689--699.
\bibitem{AB:HopfAlgSS}
\bysame,
On the cohomology of some Hopf algebroids and Hattori-Stong
theorems,
Homology, Homotopy and Applications {\bf2} (2000), 29--40.
\bibitem{AB:EnEn}
\bysame,
$I_n$-local Johnson-Wilson spectra and their Hopf algebroids,
Documenta Math. {\bf5} (2000), 351--364.
\bibitem{BCRS}
A.~Baker, F.~Clarke, N.~Ray \& L.~Schwartz,
On the Kummer congruences and the stable homotopy of $BU$,
Trans. Amer. Math. Soc. {\bf316} (1989), 385--432.
\bibitem{AB&UW}
A.~Baker \& U.~W\"urgler,
Liftings of formal group laws and the Artinian completion of $v_n^{-1}BP$,
Proc. Camb. Phil. Soc. {\bf106} (1989), 511--30.
\bibitem{BR}
M.~Basterra \& B.~Richter,
(Co-)homology theories for commutative ($S$-)algebras,
`Structured Ring Spectra', A.~Baker \& B.~Richter (eds),  
London Math. Soc. Lecture Notes Series {\bf 315} (2004), 115--131.
\bibitem{Bourbaki:CommAlg}
N.~Bourbaki,
Elements of Mathematics: Commutative Algebra, Chapters 1--7,
reprint of the 1989 English translation,
Springer-Verlag (1998).
\bibitem{CCW}
F.~Clarke, M.~D.~Crossley \& S.~Whitehouse,
Bases for cooperations in $K$-theory,
$K$-theory {\bf23} (2001), 237--250.
\bibitem{EKMM}
A.~Elmendorf, I.~Kriz, M.~Mandell \& J.~P.~May,
Rings, modules, and algebras in stable homotopy theory,
Mathematical Surveys and Monographs {\bf47} (1999).
\bibitem{GH}
P.~G.~Goerss \& M.~J.~Hopkins,
Moduli spaces of commutative ring spectra,
`Structured Ring Spectra',  A.~Baker \& B.~Richter (eds),  
London Math. Soc. Lecture Notes Series {\bf315} (2004),
151--200.
\bibitem{HS}
M.~Hovey \& N.~P.~Strickland,
Morava $K$-theories and localisations,
Memoirs of the Amer. Math. Soc. {\bf139} no.~666 (1999).
\bibitem{Jensen:Paper}
C.~U.~Jensen,
On the vanishing of $\varprojlim ^{(i)}$,
J. Algebra {\bf15} (1970), 151--166.
\bibitem{Jensen:LNM}
\bysame,
Les foncteurs d\'eriv\'es de $\varprojlim$ et leurs applications
en th\'eorie des modules,
Lecture Notes in Math. {\bf 254} (1972).
\bibitem{Joachim}
M.~Joachim,
A symmetric ring spectrum representing $KO$-theory,
Topology {\bf40} (2001), 299--308.
\bibitem{Johnson:K*K-Basis}
K.~Johnson,
The action of the stable operations of complex $K$-theory
on coefficient groups,
Illinois J. Math. {\bf28} (1984), 57--63.
\bibitem{LRS}
P.~S.~Landweber, D.~C.~Ravenel \& R.~E.~Stong,
Periodic cohomology theories defined by elliptic curves, 
 The \v{C}ech centennial (Boston, MA, 1993),  M.~Cenkl (ed.) et al.,
Contemp. Math. {\bf181} (1995), 317--337.
\bibitem{Loday}
J.-L.~Loday,
Cyclic Homology, second edition,
Grundlehren der Mathematischen Wissenschaften
{\bf301} Springer-Verlag (1998).
\bibitem{May:EinftyRings}
J.~P.~May,
$E_\infty$ ring spaces and $E_\infty$ ring spectra,
with contributions by F.~Quinn, N.~Ray and J.~Tornehave,
Lecture Notes in Mathematics {\bf577} (1977).
\bibitem{McC&S}
J.~E.~McClure \& R.~E.~Staffeldt,
On the topological Hochschild homology of $bu$,~I,
Amer. J. of Math. {\bf115} (1993), 1--45.
\bibitem{Rezk}
C.~Rezk,
Notes on the Hopkins-Miller theorem, 
Homotopy theory via algebraic geometry
and group representations. Proceedings of a conference on homotopy
theory, Evanston, IL, USA, March 23--27, 1997, M.~Mahowald (ed.) et al., 
Contemp. Math. {\bf 220} (1998), 313--366.
\bibitem{BR&Ro}
B.~Richter \& A.~Robinson,
Gamma-homology of group algebras and of polynomial algebras,
in the Proceedings of the 2002 Northwestern Conference
on Algebraic Topology, P.~Goerss \& S.~Priddy (eds),
Contemp. Math. {\bf346} (2004), 453--461.
\bibitem{Ro:Ainfty}
A.~Robinson,
Obstruction theory and the strict associativity of Morava
$K$-theories,
in `Advances in Homotopy Theory', B.~Steer, W.~Sutherland (eds), 
London Mathematical Society Lecture Note Series {\bf139}
(1989), 143--52.
\bibitem{Ro:Einfty}
\bysame,
Gamma homology, Lie representations and $E_\infty$ multiplications,
Invent. Math. {\bf152} (2003), 331--348
\bibitem{Ro:unit}
\bysame,
Classical obstructions and $S$-algebras,
`Structured Ring Spectra',  A.~Baker \& B.~Richter (eds), 
London Math. Soc. Lecture Notes Series {\bf315} (2004), 133--49.
\bibitem{RoWh}
A.~Robinson \& S.~Whitehouse,
Operads and $\Gamma$-homology of commutative rings,
Math. Proc. Cambridge Philos. Soc. {\bf 132} (2002), 197--234.
\bibitem{SVW}
R.~Schw\"anzl, R.~M.~Vogt \& F.~Waldhausen,
Adjoining roots of unity to $E_\infty$ ring spectra in good
cases -- a remark,  Homotopy invariant
algebraic structures. A conference in honor of J. Michael
Boardman. AMS  special session on homotopy theory, Baltimore, MD, USA,
January 7-10, 1998, J.-P.~Meyer (ed.) et al.,
Contemp. Math. {\bf 239} (1999), 245--249.
\bibitem{Shatz}
S.~Shatz,
Profinite Groups, Arithmetic and Geometry,
Annals of Mathematics Studies {\bf 67} (1972),
Princeton University Press.
\bibitem{Weibel}
C.~A.~Weibel,
An Introduction to Homological Algebra,
Cambridge University Press (1994).
\bibitem{Zelinsky}
D. Zelinsky,
Linearly compact modules and rings,
Amer. J. Math. {\bf75} (1953), 79--90.
\end{thebibliography}
\end{document}